\documentclass[preprint,12pt,authoryear]{elsarticle}
\usepackage{newtxtext,newtxmath} 
\usepackage{hyperref}
\usepackage{float}
\usepackage{subfigure}
\usepackage{verbatim}
\newtheorem{theorem}{Theorem}
\newtheorem{prop}{Proposition}
\newtheorem{remark}{Remark}
\newtheorem{proof}{Proof}

\journal{}

\begin{document}
\begin{frontmatter}


\title{Experimental designs for controlling the correlation of estimators in two parameter models}


\author[inst1,inst2]{Edgar Benitez}

\affiliation[inst1]{organization={Datai, Institute for Data Science and
Artificial Intelligence, University of Navarra},
            addressline={Campus Universitario}, 
            city={Pamplona},
            postcode={31009}, 
            state={Navarra},
            country={Spain}}
\affiliation[inst2]{organization={Tecnun Escuela de Ingeniería, University of Navarra},
            addressline={Campus Universitario}, 
            city={Pamplona},
            postcode={31009}, 
            state={Navarra},
            country={Spain}}
\author[inst1,inst2]{Jesús López-Fidalgo}

\begin{abstract}
The state of the art related to parameter correlation in two-parameter models has been reviewed in this paper. The apparent contradictions between the different authors regarding the ability of  D--optimality to simultaneously reduce the correlation and the area of the confidence ellipse in two-parameter models were analyzed. Two main approaches were found: 1) those who consider that the optimality criteria simultaneously control the precision and correlation of the parameter estimators; and 2) those that consider a combination of criteria to achieve the same objective. An analytical criterion combining in its structure both the optimality of the precision of the estimators of the parameters and the reduction of the correlation between their estimators is provided. The criterion was tested both in a simple linear regression model, considering all possible design spaces, and in a non-linear model with  strong correlation of the estimators of the parameters (Michaelis--Menten) to show its performance. This criterion showed a superior behavior to all the strategies and criteria to control at the same time the precision and the correlation.
\end{abstract}


\end{frontmatter}



\section{Introduction}

\subsection{Optimal estimators for two parameter models}
Let the model
\begin{equation}
y=f^{T}(x)\pmb{\theta}+\varepsilon,\; x\in\chi,
\end{equation}
where $\pmb{\theta}$ is the vector of parameters to be estimated, \textit{f(x)} the vector of regressors, $\chi$ the design space and the error $\varepsilon$ is assumed to be Normal with constant variance. The estimation of the parameters is made performing experiments under the conditions $x_{1},\dots,x_{n}$ and obtaining the corresponding responses $y_{1},\dots,y_{n}$. The collection $x_{1},\dots,x_{n}$ is called an {\it exact experimental design of size} \textit{n}. Some of these points may be repeated and a probability distribution can be defined by assigning to each point the proportion of times it appears in the design as a probability measure. Thus, the more general concept of {\it approximate design} can be given as any measure of probability, $\xi$.

The Fisher Information Matrix (FIM) of a design $\xi$ is
\begin{equation}
M(\xi)=\int{f(x)^{T}f(x) \xi(dx)}.
\end{equation}

Its inverse is proportional to the covariance matrix of the estimators of the parameters,
\begin{equation*}
\Sigma =\frac{\sigma^{2}}{n}M^{-1}\left ( \xi \right ),
\end{equation*}
which for two parameters will be 
\begin{equation*}
\Sigma =\begin{pmatrix}
\sigma_{1}^{2} & \sigma_{12} \\
\sigma_{12} & \sigma_{2}^{2} 
\end{pmatrix}.
\end{equation*}
From the FIM, scalars are used for representing different mathematical properties of the estimated model. Properties of statistical interest relates to precision or correlation of the estimators of the parameters, which can be optimized for particular design conditions $x_{i},\dots,x_{n}$. These features of the FIM are called {\it optimality criteria}, $\Phi[M(\xi)]$ or just $\Phi(\xi)$ for simplicity. This function is assumed non-increasing in the Loewner sense, that is, if $M(\xi)-M(\xi^prime)$ is non--definite negative then $\Phi[M(\xi)]\leq \Phi[M(\xi^prime)]$. A $\Phi$--optimal design is a design minimizing this function. Typically $\Phi$ is convex and frequently differentiable. Then, with a few other conditions an {\it equivalence theorem} can be proved in order to check whether a particular design, $\xi^{\star}$, is optimal,
\begin{equation*}
\partial \Phi[M(\xi^{\star}),M(\xi_x)] \geq 0,\, x\in \chi,
\end{equation*}
where $\xi_x$ stands for a one--point design giving the whole mass to $x$. The directional derivative is defined as
\begin{equation*}
\partial \Phi[M(\xi),M(\xi^{\prime})]=\lim_{\alpha \rightarrow 0^+} \frac{\Phi[(1-\alpha)M(\xi)+\alpha M(\xi^{\prime})]-\Phi[M(\xi)]}{\alpha}.
\end{equation*}

This theorem can also be used for generating algorithms for computing optimal designs. It is only valid for approximate designs. If the criterion is inverse-positively homogeneous in the sense that $\Phi[\lambda M(\xi)]=\frac{1}{\lambda}\Phi[M(\xi)]$ for any $\lambda$ then the efficiency of a design $\xi$,
\begin{equation*}
\text{Eff}_{\Phi} (\xi)=\frac{\Phi[M(\xi^{\star})]}{\Phi[M(\xi)]},
\end{equation*}
is a number between 0 and 1, usually multiplied by 100 and referred in terms of percentages,  with a practical interpretation. For instance, an efficiency of 70\% means 30\% of the experiments (observations) can be saved using the optimal design $\xi^{\star}$ instead of $\xi$.

\begin{figure}
\includegraphics[trim={1mm 1mm 1mm 1mm},clip,width=12cm]{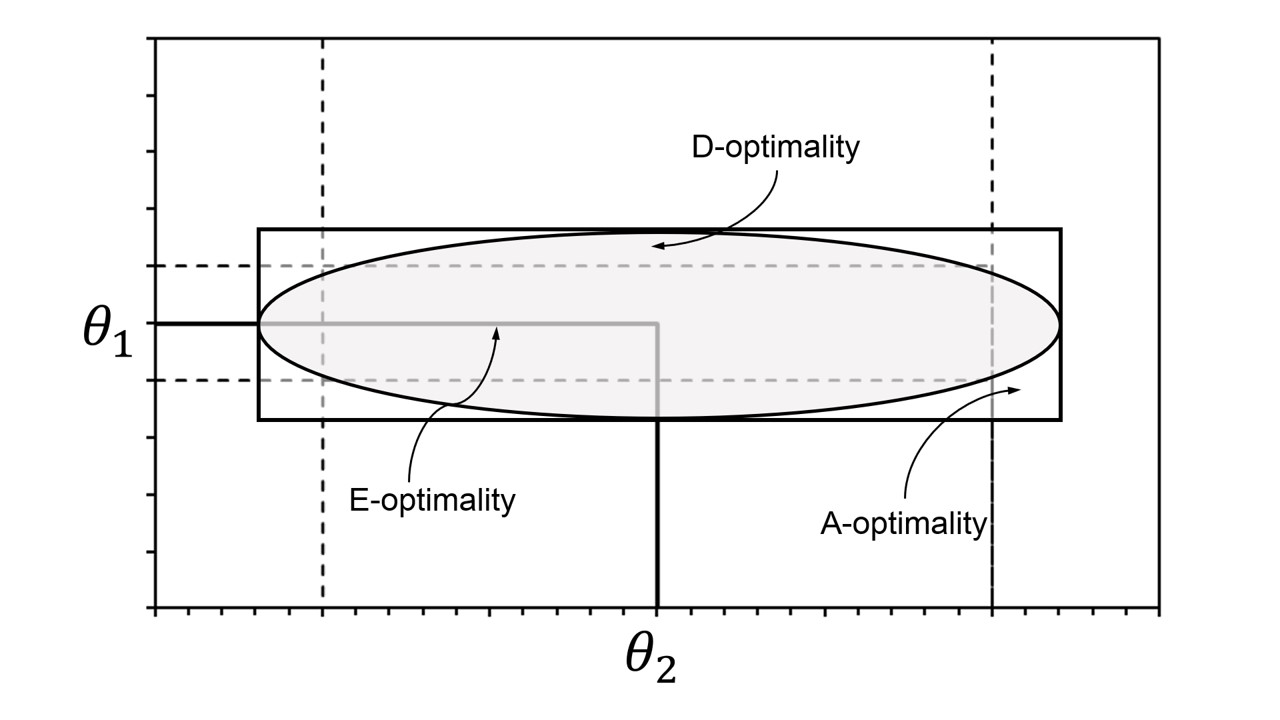}
\caption{Confidence ellipse of a two-parameter model}
\label{f01e2p}
\end{figure}

Some strategies can be evaluated geometrically in relation to the confidence ellipsoid of the parameters. Figure \ref{f01e2p} shows the confidence ellipse for a model with two parameters in the particular case where the estimators are uncorrelated. For correlated estimators the ellipse would be in a inclined position. Three typical criteria such as D--, A-- and E--optimality, look for ``minimizing'' in different ways the geometry of the ellipse. In particular, D--optimality minimizes the area of the ellipse while E--optimality minimizes the largest axis of it and A--optimality is trying to minimize the average of the axes in the uncorrelated case. The dashed lines in Figure \ref{f01e2p} show the marginal confidence intervals of each parameter.

Some of the first works on designs for correlated parameter estimators considered that a ``convenient approach" for minimizing the correlation is the use of the determinant of the variance-covariance matrix \cite{duggleby, box_lucas}. It seems that some authors from this took for granted, wrongly, that D--optimality minimizes the correlation between the estimators for models with two parameters. These are some references,

\begin{quote}
    {\it ``Box and Hunter \cite{boxhunter} suggested that for a two-parameter model the volume criterion (D-optimal design) minimizes not only the size of the uncertainty region associated with the parameter estimates but also the correlation between them.''} \cite{Franceschini}
    \end{quote}

\begin{quote}
{\it ``The D-optimal design can both reduce the parameter correlation and improve the parameter precision for two-parameter system, but the approach cannot be extended to model systems with three or more parameters.''} \cite{stateart}
\end{quote}

\begin{quote}
{\it ``One of the most widely used sequential experimental design criteria for nonlinear modelling is minimization of the determinant of the covariance matrix for the parameter estimates, $|V|$, recommended by Box and Lucas and many others have recently reviewed this criterion and its extensive applications. Although it is recognized that sets of parameters estimate with the same determinant $|V|$ can have very different correlation structures, Box and Hunter suggested that for a two-parameter model the criterion minimizes not only the size of the region of uncertainty associated with the parameter estimates but also the correlation between them However, the determinant design criterion does not necessarily minimize correlations among parameter estimates in models with three or more parameters.''} \cite{pritchard_bacon}
\end{quote}

Actually, Box and Hunter \cite{boxhunter} established a procedure to compute optimal designs sequentially for non--linear models. They were aware of the correlation between estimators for some models, especially in catalytic reaction kinetics. Thus, they proposed to check the correlation at some steps of the sequential procedure adding points to reduce the correlation,

\begin{quote}
{\it ``(...) there may be a special reason for wanting to minimize the correlation between a certain pair of estimates. In such situations, one could proceed by calculating the variance and covariance terms which are of special interest in addition to using the overall criterion.''} \cite{boxhunter} 
\end{quote}

However, this claim has been dismissed by different studies \cite{bhonsale, Telen}, where none of the proposed criteria meets the purposed condition (efficient parameter estimation and what with call hereafter {\it decorrelation}). As a matter of fact, minimizing the determinant of a $2 \times 2$ covariance matrix means minimizing the product of the variances while maximizing the covariance,

\begin{equation}
\label{eq:01_dbl}
\left | \Sigma \right |=\sigma_{1}^{2}\sigma_{2}^{2}-\sigma_{12}^{2}.
\end{equation}

Nevertheless, minimizing the correlation instead of the variance should be the final objective. Both optimization problems are quite different as it will be seen later. Additionally, minimizing the correlation takes frequently to singular designs, which is a strong reason for using it beside other criteria. Compound criteria are then a typical choice, but the criterion of the product of the variances takes care at the same time of minimizing the square of the correlation and the determinant of the inverse of the information matrix.

\subsection{Why is it important to consider the correlation of the estimators?}
In the cases of models with a single parameter there is just one universal optimality criterion. Otherwise, for models with several parameters, in addition to the precision of the estimators of the parameters, the experimental design is also sought to control the level of correlation between them (decorrelation). The study of the decorrelation of the estimators of the parameters, in general, shows an extensive work in the literature with a number of reviews and proposals. 

Minimally correlated estimators of the parameters are important for a number of reasons. On the one hand, correlated estimators imply linear dependence of the parameter estimators and the FIM becomes near singular, which makes the calculation of the inverse computationally problematic. This may lead to an important computational problem of the maximum likelihood estimators (MLE), which for linear models are the least squared estimators (LSE). Additionally, this introduces a problem of parameter identification with a confounding overlapping of the actual effects of the regressors \cite{mclean}. Thus, inaccurate estimators with unreliable statistical tests come straightforward. Ultimately, this is reflected in the low t value and therefore leads to the non--rejection of the null hypothesis for null parameters \cite{stateart}.

Some authors consider that the parameters are not estimable for a cut-off point in correlation of $\pm$0.99 \cite{valcor}. However, it is not necessary to have high correlation between parameter estimators to affect the design efficiency. In particular, the ``offset" effect may appear when a change in a parameter is compensated by changes in the other(s) leading to different values of the parameter estimators with similar values of the objective function \cite{corr}.

\begin{figure}
\subfigure[No correlated]{\includegraphics[width=65mm]{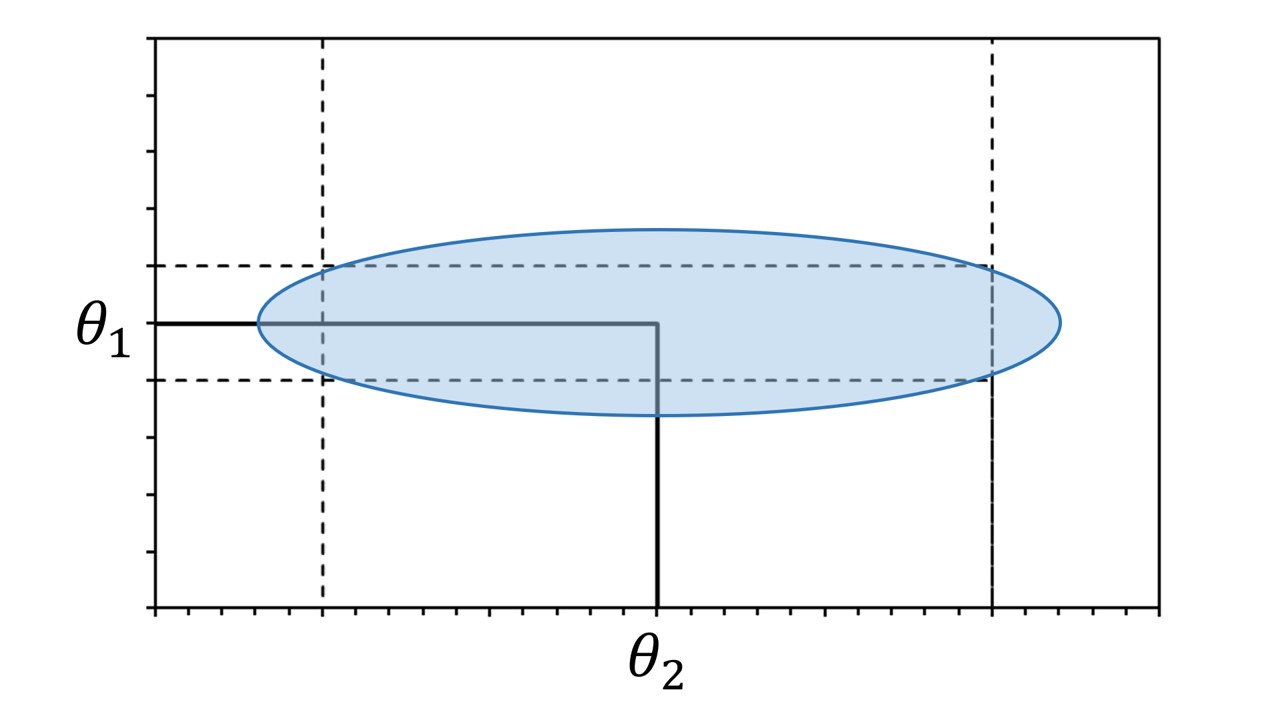}}
\subfigure[Correlated]{\includegraphics[width=65mm]{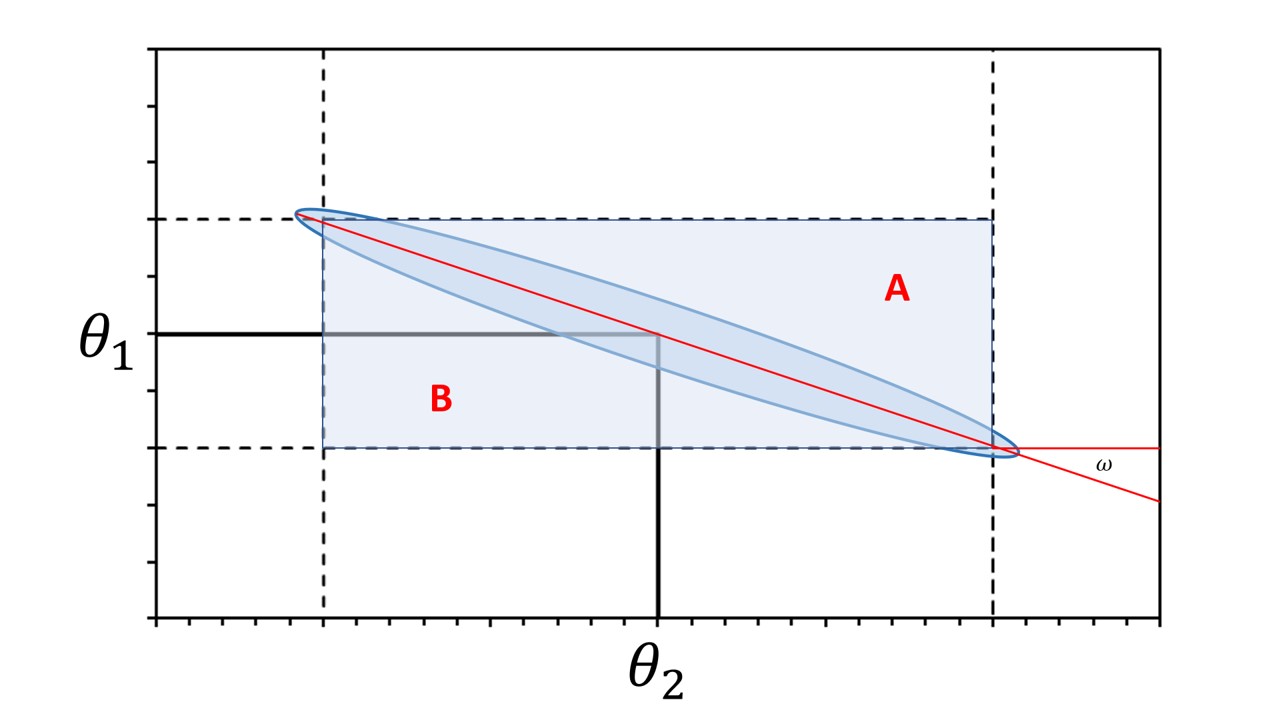}}
\caption{Ellipsoid of the confidence intervals of the estimators of a model with two parameters. Dotted lines represent the marginal confidence intervals for each parameter}
\label{f03ecn}
\end{figure}

In summary, the associated problems can be classified into three types: 1) Estimability: the estimators of the parameters are confounded; 2) some t-values might be low just due to the high correlation between the parameter estimators \cite{wang, maheshwari}; and 3) for the same values of the objective function there are different combinations of the estimators of the parameters. Figure \ref{f03ecn} shows confidence regions for uncorrelated (a) and correlated estimators (b). The shaded rectangles are formed by the marginal confidence intervals of the parameters. Thus, the area of this rectangle is much larger for the correlated estimators. For (a) the area of the confidence rectangle coincides almost with the area of the ellipsoid. On the contrary, for (b) the area of the rectangle (regions A and B) is much larger than the area of the ellipse. Therefore, in the case of uncorrelated parameters, the joint confidence interval is a good approximation to the real confidence ellipsoid of the estimation of the parameters, which does not occur with correlated parameter estimators. This means it is difficult to know how well each parameter will be estimated individually. For two parameters a measure of this degree of uncertainty can be quantified by the relationship between the two axes of the ellipsoid, $\theta_{1}\cos \omega+\theta_{2}\sin \omega$ for the largest and $-\theta_{1}\sin\omega+\theta_{2}\cos\omega$ for the shortest, where $\omega$ is the angle between the largest axis of the ellipsoid and the x-axis.

It can be assessed that the correlation of parameter estimators is an undesirable characteristic in a model, both from the phenomenological point of view and from its mathematical and statistical properties mentioned above. From the perspective of the phenomenon being measured, it is clear that models with multiple parameters can be useful, since each one represents natural processes that can be described in terms of a specific theoretical explanation of interest to the researcher. The fact that some estimators of a model show consistently correlation with each other, would indicate that this assumption is not entirely true, and in that sense phenomena that are considered independent could merge under the estimation of the corresponding parameters. 
\subsection{How to deal with correlation problems?}
Solutions for the estimators correlation issue are model reparameterization \cite{agarwal}, variable separation, independent estimation of the parameters, or using appropriate experimental designs. In some studies with models with two highly correlated parameters, such as tumor growth models, this linear relationship is previously estimated. Then, one parameter is worked out in that expression, say, for instance $\theta_2= a_0 +b_0 \theta_1$, and plugged into the model with a reduction of the dimensionality of the parameters. Some authors showed that these simplified versions of the models were just as efficient in prognosis as the versions with the two parameters \cite{gomred} for high correlation. Nevertheless, to date, the most satisfactory strategy is the use of good experimental designs \cite{box_lucas}. 

{In general, optimality criteria have a wide range of effects on the reduction of the correlation of the parameter estimators. On the one hand, A-optimality does not take into account the entire covariance matrix, just the trace, that is the sum of variances. Thus, it seems insensitive to correlation \cite{Franceschini}}. On the other hand, E-optimality takes the longest axis of the confidence ellipsoid as the minimization criterion. However, more efficient modifications of this criterion take into account not only this axis, but also its relationship with the minimum axis (such as the Modified E- or EM--optimality), for which the theoretical optimum value is one, which is equivalent to minimizing the condition number (ratio between the largest and the minor axes of the ellipsoid). Far from being accepted, this criterion falls back into the problem of the previous ones, since it determines the search for spheroidal shapes, which tend to generate large volumes in confidence ellipsoids.

As mentioned in the Introduction, simultaneous control properties have been assigned to the D-optimality, both in the precision of the estimation and in the correlation of the estimators of the parameters. However, this must be clarified, since it is not strictly true. 

\subsection{Compound and multi--criteria selection}
What can be deduced from the previously mentioned works is that the strategy of obtaining an indicator that achieves a simultaneous optimization of both, the estimation of the parameters and their correlations has not been found successfully. Therefore, the attempt to create multi--criteria indicators to solve the problem has generated wide interest \cite{bhonsale,Telen}. Under this approach, two types of optimality can be distinguished: 1) Compound optimality criteria, which combine properties of more than one criterion, for example, a weighted sum of two criteria such as $\phi_{D}$ and $\phi_{EM}$ \cite{bhonsale}; $\phi_{A}$ and $\phi_{EM}$ \cite{wang}; or $\phi_{c}$ with restrictions on eigenvalues \cite{wang, Franceschini} and 2) Select from a set of optimal solutions using Pareto fronts \cite{Telen,bhonsale,das}.
Some authors have highlighted that taking into account A--optimality only contains the information of the parameters \textit{per se}, then this criterion seems very convenient to be combined with another criterion that gives weight to the correlation of parameters \cite{Franceschini}.

\section{A proposal for decorrelation of the estimators of the parameters}
In this section the criterion of the product of the variances of the estimators for models with two parameters is going to be related to the correlation of the estimators and the determinant of the information matrix.

For a model with $m$ parameters the criterion function for D--optimality is defined as
\begin{gather*}
\Phi_{D}[M(\xi)]=\left | M^{-1/m}(\xi) \right |
\end{gather*}
if the information matrix is non-singular. 

Dette (1997) introduced R-optimality, the product of the variances of the MLE's (the diagonal of the inverse of the FIM) as a criterion minimizing the area of a confidence hyper-rectangle of the parameters based on Bonferroni t-intervals \cite{dette_1}. He gave the criterion of \textit{R--optimallity} in a context of standardized criteria to make the variances comparable using efficiencies for estimating each parameter. Neither D-- nor R--optimality were able to be standardized in this way. He gave this criterion just as an alternative to D--optimality. As far as authors' knowledge this criterion has not been used in the posterior literature. It is defined as the product of the variances of the estimators,
\begin{gather*}
\Phi_{R}(\xi)=\prod _{i=1}^{m}\left[ M(\xi) \right ]_{ii}^{-1/m}.
\end{gather*}

In his work Dette provided some results for a class of criteria including standardized A- and MV-optimality as well as R-optimality. 

For two parameters, $m=2$, these criteria are 
\begin{gather*}
\Phi_{D}[M(\xi)]=\left | M^{-1/2}(\xi) \right |=\left[\frac{n}{\sigma^2}(\sigma_{1}^2\;\sigma^2_{2}-\sigma_{12}^{2})\right]^{1/2}\\
\Phi^2_{R}[M(\xi)]=
\left(\frac{n}{\sigma^2} \sigma_{1}^2\;\sigma_{2}^2\right)^{1/2},
\end{gather*}
 and the criterion of the square of the correlation will be called $r^2$--optimality,
\begin{gather*}
\Phi_{r^2}[M(\xi)]=\frac{\sigma_{12}^{2}}{\sigma_{1}^2\;\sigma_{2}^2}.
\end{gather*}
From which this relationship arises directly:
\begin{equation}
\label{eq:03_cdr}
\Phi^2_{R}=\frac{\Phi^2_{D}}{1-\Phi_{r^{2}}},
\end{equation}

Therefore, the R criterion controls simultaneously both the precision in the estimation of the parameters and the correlation between them.

\subsection{Convexity and differentiability of the R--optimal criterion}
In order to apply the equivalence theorem to this criterion the convexity must be proved. Then the differentiability is proved and the directional derivative computed in this section.

\begin{theorem}
Let R--optimality for two parameters be defined as the inverse of the product of the variances of the estimators for a non--singular design, $\xi$,
$$\Phi_R^2(\xi)=\{M^{-1}(\xi)\}_{11}\{M^{-1}(\xi)\}_{22}.$$
\begin{enumerate}
 \item[i)] This function is convex.
  \item[ii)] This function is differentiable and the gradient is 
  \begin{eqnarray*}
\bigtriangledown \Phi_R(\xi)&=&\frac{-1}{2\Phi_R(\xi)}\left[\frac{M^{-1}(\xi)}{\left | M^{-1}(\xi) \right |}+M^{-1}(\xi)c_2c_1^T M^{-1}(\xi)\right].  
  \end{eqnarray*}
  
Therefore the directional derivative is
\small
  \begin{eqnarray*}
\frac{1}{2\Phi_R(\xi)}[2\left | M(\xi) \right |&+&c_1^TM^{-1}(\xi)c_2)-f^T(x)M^{-1}(\xi)f(x)\left | M(\xi) \right |\\&-&f^T(x)M^{-1}(\xi)c_2c_1^TM^{-1}(\xi)f(x)].
  \end{eqnarray*} 
  \end{enumerate}
\end{theorem}

\begin{proof}
\begin{enumerate}
\item[(i)] Let $\alpha \in (0,1)$ and $\xi_1$, $\xi_2$ be two non--singular designs. For simplicity of notation let us call 
$$M^{-1}(\xi_1)=\left(\begin{array}{cc}
 u_1 & u_{12} \\
 u_{12} & u_2
\end{array}\right), \quad
M^{-1}(\xi_2)=\left(\begin{array}{cc}
 v_1 & v_{12} \\
 v_{12} & v_2
\end{array}\right).$$
It is known that the c--criterion function $\Phi_c(\xi)=c^T M^{-1}(\xi)c$, is convex for non-singular designs. In particular it is convex for $c_1=(1,0)^T$ and $c_2=(0,1)^T$: $\Phi_{c_1}(\xi)=c_1^T M^{-1}(\xi)c_1=\{M^{-1}(\xi)\}_{11}$ and $\Phi_{c_2}(\xi)=c_2^T M^{-1}(\xi)c_2=\{M^{-1}(\xi)\}_{22}$.
\begin{enumerate}
\item[a)] The result will be proved first for the case $u_1 \geq v_1$ and $u_2 \geq v_2$ (or the symmetrical case). Using the convexity of c-optimality,
\small
\begin{eqnarray*}
\Phi_R^2[(1-\alpha)\xi_1+\alpha \xi_2]&=&\{M^{-1}((1-\alpha)\xi_1+\alpha \xi_2)\}_{11}\{M^{-1}((1-\alpha)\xi_1+\alpha \xi_2)\}_{22}\\ &\leq& [(1-\alpha)u_1+\alpha v_1][(1-\alpha)u_2+\alpha v_2]\leq (1-\alpha)u_1 u_2+\alpha v_1 v_2.  
  \end{eqnarray*}

  The last inequality arises by taking into account that $(u_2-v_2)(u_1-v_1)\geq 0$ and therefore $u_1 (v_2-u_2)+v_1(u_2-v_2)\leq 0$. Then $[(1-\alpha)^2-(1-\alpha)]u_1u_2+[\alpha^2-\alpha]v_1v_2+\alpha (1-\alpha)(u_1v-2+u_2v_1)\leq 0$ and the inequality is proved.

 \item[b)] Let assume now that  $u_1 \leq v_1$ and $u_2 \geq v_2$ (or the symmetrical case). Then
    \begin{eqnarray*}
\{M^{-1}((1-\alpha)\xi_1+\alpha \xi_2)\}_{11}&\leq & \{M^{-1}(\xi_1)\}_{11}=u_1,\\
\{M^{-1}((1-\alpha)\xi_1+\alpha \xi_2)\}_{22}&\leq & \{M^{-1}(\xi_2)\}_{22}=v_2.
  \end{eqnarray*}
Thus,
    \begin{eqnarray*}
\Phi_R^2[(1-\alpha)\xi_1+\alpha \xi_2]&\leq& u1v2 \leq (1-\alpha)u_1u_2+\alpha v_1 v_2.
  \end{eqnarray*}
  Last inequality comes from the fact that 
      \begin{eqnarray*}
(1-\alpha)(u_1u_2-u_1v_2)+\alpha (v_1 v_2- u_1 v_2) = (1-\alpha)u_1(u_2-v_2)+\alpha v_2 (v_1-u_1) \geq 0.  
  \end{eqnarray*}
\end{enumerate}
  
 \item[(ii)] This function can be expressed as the sum of two differentiable functions,
$$\Phi_R^2(\xi)=\Phi_D^2(\xi)+c_1^TM^{-1}(\xi)c_2.$$
A result is that $\bigtriangledown \text{tr}AM^{-1}B =-M^{-1}A^TB^T M^{-1}$.
Then the gradient is,
  \begin{eqnarray*}
\bigtriangledown \Phi_R^2(\xi)&=&\bigtriangledown \Phi_D^2(\xi)+ \bigtriangledown c_1^TM^{-1}(\xi)c_2\\
&=& \frac{-M^{-1}(\xi)}{\left | M^{-1}(\xi) \right |}-M^{-1}(\xi)c_2c_1^T M^{-1}(\xi).  
  \end{eqnarray*}
  Therefore,
  \begin{eqnarray*}
\bigtriangledown \Phi_R(\xi)&=&\frac{-1}{2\Phi_R(\xi)}\left[\frac{M^{-1}(\xi)}{\left | M^{-1}(\xi) \right |}+M^{-1}(\xi)c_2c_1^T M^{-1}(\xi)\right].  
  \end{eqnarray*}

  The directional derivative in the direction of a one-point design is then
  \begin{eqnarray*}
\partial \Phi_R(\xi,\xi_x)&=&\frac{1}{2\Phi_R(\xi)\left | M^{-1}(\xi) \right |}+\frac{\text{tr}(M^{-1}(\xi)c_2c_1^T)}{2\Phi_R(\xi)} \\ &&-\frac{\frac{f^T(x)(M^{-1}(\xi)f(x)}{\left | M^{-1}(\xi) \right |}+f^T(x)M^{-1}(\xi)c_2c_1^TM^{-1}(\xi)f(x)}{2\Phi_R(\xi)}.
  \end{eqnarray*} 
\end{enumerate}
\end{proof}

\section{A complete study of correlation for Simple linear regression}
Let 
  \begin{eqnarray*}
y&=&\theta_1+\theta_2 x + \varepsilon,\quad \varepsilon \sim \mathcal{N} (0,\sigma^2),\quad x\in \chi =[a,b],
  \end{eqnarray*}
with independent observations. 

The Information Matrix for  an approximate design with finite support,
  \begin{eqnarray*}
\xi&=&\left\{\begin{array}{ccc}
 x_1 & \cdots & x_k \\
 p_1 & \cdots & p_k 
\end{array}\right\},
  \end{eqnarray*}
is
   \begin{eqnarray*}
M(\xi)&=&\sum_i p_i \left(\begin{array}{cc}
 1 & x_i \\
 x_i & x_i^2 
\end{array}\right) = \left(\begin{array}{cc}
 1 & \overline{x} \\
 \overline{x} & \overline{x^2} 
\end{array}\right),
  \end{eqnarray*}
  where $\overline{x^r}=\sum_i p_i x_i^r$ for any integer $r$. The inverse is then
   \begin{eqnarray*}
M^{-1}(\xi)&=&\frac{1}{S_x^2} \left(\begin{array}{cc}
 \overline{x^2} & -\overline{x} \\
 -\overline{x} & 1 
\end{array}\right),
  \end{eqnarray*} 
  where $S_x^2=\sum_i p_i (x_i-\overline{x})^2$.

We want to analyse exhaustively the correlation between the MLE's of $\theta_1$ and $\theta_2$ for this model.
  \subsection{D-optimal design}

It is well known that the D--optimal design for this model is
   \begin{eqnarray*}
\xi_D&=&\left\{\begin{array}{cc}
 a & b \\
 1/2 & 1/2 
\end{array}\right\}. 
  \end{eqnarray*}

   \subsection{R-optimal design}
  The directional derivative of the R--optimality criterion  for simple linear regression is a second order polynomial, which means the R--optimal design always has two points and they should be the extreme values of the interval. The computation of the optimal weights 
  come from solving the equation where the derivative with respect to the second weight vanishes,
   \begin{eqnarray*}
\xi_R=\left\{\begin{array}{cc}
 a & b \\
 1-p_R & p_R 
\end{array}\right\}
  \end{eqnarray*}
 where 
   \begin{eqnarray*}
p_R= \left\{\begin{array}{cc} 
\frac{4 a^2}{5 a^2+A-b^2} & \text{ if } a, b \neq 0,\\
1/3 & \text{ if } a= 0,\\
2/3 & \text{ if } b= 0,
\end{array}\right.
   \end{eqnarray*}
and $A=\sqrt{a^4+14 a^2 b^2+b^4}$, for simplicity of notation. For $b=-a$, $p_R=1/2$.

  \subsection{Optimal correlation design}
 The criterion of the square correlation (determination coefficient) is defined for non-singular designs as 
 $$\Phi_{r^2}(\xi)=\frac{\{M^{-1}(\xi)\}_{12}^2}{\{M^{-1}(\xi)\}_{11}\{M^{-1}(\xi)\}_{22}}.$$

 This criterion is not convex and therefore the equivalence theorem cannot be applied. Caratheodory's theorem is still applicable, thus there is an optimal design with no more than 4 points. 

 For a 2--point design,
  \begin{eqnarray*}
\xi&=&\left\{\begin{array}{cc}
 x_1 & x_2 \\
 1-p & p 
\end{array}\right\},
  \end{eqnarray*}
the square of the correlation for simple linear regression is 
 \begin{eqnarray}
\frac{[(1-p) x_1+p x_2]^2}{(1-p) x_1^2-p x_2^2}. \label{eq.corr2}
  \end{eqnarray}

 If $0\leq a<b$ or $a<b\leq 0$ then the 2--point design with the smallest correlation is
   \begin{eqnarray*}
\xi_{r^2}=\left\{\begin{array}{cc}
 a & b \\
 \frac{|b|}{|a|+|b|} & \frac{|a|}{|a|+|b|} 
\end{array}\right\},
  \end{eqnarray*}
Notice that for $a=0$ or $b=0$ this is a one--point design at zero and then the slope of the regression line cannot be estimated and therefore the correlation cannot be computed. If  one of the points in  (\ref{eq.corr2}) is zero then the correlation is the weight of the second point, no matter which point it is. The message is that taking a two-point design with zero and other point will give a square correlation equal to the weight of this second point. Fixing this weight the optimal point can be optimized according to D-- or R--optimality, which in both cases takes to the other extreme value of the interval.

If $a<0<b$ then any design $\xi$ such that $\overline{x}=\sum_x \xi(x) x =0$ produces zero correlation and then they are $r^2$--optimal designs. Let $\Xi_R=\{\xi \,|\, \overline{x}=\sum_x \xi(x) x =0 \}$, which is convex and then the equivalence theorem can be applied in it. Any set of points, not all with the same sign, may produce a design in $\Xi_R$ with appropriate weights. Moreover, in this space both D-- and R--optimality functions are the same and therefore the restricted optimal designs. The R-- and D--optimality criteria restricted to $\Xi_R$ are equal, $\xi_{D}=\xi_{R}$.

\subsection{D--efficiencies}
  \begin{prop}
 \begin{enumerate}
  \item[i)] The D-efficiency of $R$--optimal designs is \end{enumerate} 
  \begin{eqnarray*}
  \text{Eff}_D (\xi_R)&=& 4 |a| \frac{\sqrt{a^2+A-b^2}}{5 a^2+A-b^2} \geq \frac{2\sqrt{2}}{3}=0.943,
   \end{eqnarray*}
 where the minimum is reached at $a=-1, b=0$ or $a=0, b=1$, the efficiency is 1 for $b=-a$ and also when $b \rightarrow a$. 
  \item[ii)] The D-efficiency of $r^2$--optimal designs is 
   \begin{eqnarray*}
\text{Eff}_D (\xi_{r^2})&=&\left\{\begin{array}{cc}
 0 &\text{if } a=0, \text{ or } b=0,\\ 
 2 \frac{\sqrt{|ab|}}{|b|+|a|} & \text{otherwise}.
\end{array}\right.
  \end{eqnarray*}
Thus, $\text{Eff}_D (\xi_{r^2}) \rightarrow 1$ as $b\rightarrow a \neq 0$.
   \end{prop}

\subsection{R--efficiencies}
\begin{prop}
\begin{enumerate}
 \item[i)] The R--efficiency of a D--optimal design is
 \begin{eqnarray*} 
 \text{Eff}_R (\xi_D)&=& \frac{1}{8} \sqrt{34-\frac{a^2}{b^2}+\frac{\left(a^2+13 b^2\right) A}{b^2 \left(a^2+b^2\right)}+\frac{A-b^2}{a^2}}\geq \frac{3 \sqrt{\frac{3}{2}}}{4} = 0.919.\\
  \end{eqnarray*}
The equality is reached for $a=-1$ and $b=0$, or $a=0$ and $b=1$. Then, $\text{Eff}_R (\xi_D)\rightarrow 1$ as $b\rightarrow a \neq 0$.
  \item[ii)] The R-efficiency of $r^2$--optimal designs is 
  \begin{eqnarray*} 
 \text{Eff}_R (\xi_{r^2})&=&\left\{\begin{array}{ll}
 0 &\text{if } a=0 \text{ or } b=0,\\ 
\frac{\left(A+5 b^2-a^2\right) \sqrt{\frac{|b|}{|a|}[3 a^4 -b^4+14 a^2b^2+A \left(b^2 +3 a^2\right)]}}{8 \sqrt{2} b^2 (|a|+|b|)^2} & \text{otherwise}.
\end{array}\right.
  \end{eqnarray*}
Again $ \text{Eff}_R (\xi_{r^2})\rightarrow 1$ as $b\rightarrow a \neq 0$.
\end{enumerate}
\end{prop}

\subsection{Correlations}
Correlations are computed directly instead of $r^2$--efficiencies since this criterion is not inverse-positively homogeneous and therefore there is not a a statistical interpretation of it while the interpretation of the correlation does not need any additional justification.

\begin{prop}
\begin{enumerate}
 \item[i)] D--optimal designs produce the following correlations:
 \begin{eqnarray*}
\text{Corr} (\xi_{D})&=& -\frac{a+b}{\sqrt{2(a^2+b^2)}}\\
  \end{eqnarray*}

$ \text{Corr} (\xi_{D})\rightarrow 0$ as $b\rightarrow -a \neq 0$; $ \text{Corr} (\xi_{D})\rightarrow \pm 1$ as $b\rightarrow a \neq 0$ and $\text{Corr} (\xi_{D})\rightarrow \pm \frac{1}{\sqrt{2}}=\pm 0.707 $ as $a$ or $b \rightarrow 0$ .
  \item[ii)] R--optimal designs produce the following correlations
 \begin{eqnarray*}
\text{Corr} (\xi_R)&=& \frac{| b| \left(a^2+4 a b+A-b^2\right) \left(-5 a^2-A+b^2\right)}{\text{sgn}(a)
 \left(a^2+A-b^2\right) \sqrt{2 A^3-2 \left(a^6-33 a^4 b^2-33 a^2 b^4+b^6\right)}}\\
  \end{eqnarray*}

$ \text{Corr} (\xi_{R})\rightarrow 0$ as $b\rightarrow -a \neq 0$; $\text{Corr} (\xi_{R})\rightarrow \frac{1}{\sqrt{3}}= 0.577 $ as $b \rightarrow 0$ and $\text{Corr} (\xi_{R})\rightarrow -\frac{1}{\sqrt{3}}=-0.577 $ as $a \rightarrow 0$.
  
   \item[iii)] $r^2$--optimal designs produce the following correlations
    \begin{eqnarray*}
\text{Corr} (\xi_{r^2})&=&\left\{\begin{array}{ll}
0 & \text{if } a\leq 0 \leq b,\\ 
-\frac{2 \sqrt{a b}}{a+b} & \text{if } 0<a<b \text{ or } a<b<0
\end{array}\right.
  \end{eqnarray*}
\end{enumerate}

\end{prop}

Table \ref{tb.effs} shows cross D-- and R--efficiencies and correlations for a design space $\chi=[a,5]$ for $a =-3, -1, -1/2, -1/5, 1/5, 1/2, 1, 3$ representing most of the typical situations.

\begin{table}
\begin{center}
\resizebox{15cm}{!}{
\begin{tabular}{cccccccccc}
\hline \\ 
$a$ &$p_R$& $p_{r^2} $& $\text{Eff}_D(\xi_R)$& $\text{Eff}_D(\xi_{r^2})$&$\text{Eff}_R(\xi_D)$&$\text{Eff}_R(\xi_{r^2})$&$\text{Corr}(\xi_D)$&$\text{Corr}(\xi_R)$&$\text{Corr}(\xi_{r^2})$\\ \hline
3 & 0.439 & 0.375 & 0.992 & 0.968 & 0.986 & 0.984 & -0.970 & -0.969 & -0.968\\ 
1 & 0.356 & 0.167 & 0.958 & 0.745 & 0.934 & 0.837 & -0.832 & -0.785 & -0.745\\ 
1/2 & 0.340 & 0.0909 & 0.947 & 0.575 & 0.923 & 0.686 & -0.774 & -0.689 & -0.575\\ 
1/5 & 0.334 & 0.0385 & 0.944 & 0.385 & 0.919 & 0.481 & -0.735 & -0.623 & -0.385\\ 
-1/5 & 0.334 & 0.0385 & 0.944 & 0.385 & 0.919 & 0.481 & -0.678 & -0.531 & 0\\
-1/2 & 0.340 & 0.0909 & 0.947 & 0.575 & 0.923 & 0.686 & -0.633 & -0.465 & 0\\ 
-1 & 0.356 & 0.167 & 0.958 & 0.745 & 0.934 & 0.837 & -0.555 & -0.367 & 0\\ 
-3 & 0.439 & 0.375 & 0.992 & 0.968 & 0.986 & 0.984 & -0.243 & -0.127 & 0\\ 
-5 & 0.500 & 0.500 & 1.00 & 1.00 & 1.00 & 1.00 & 0 & 0 & 0\\ \hline 
\end{tabular} 
}
\caption{Cross D-- and R--efficiencies and correlations for simple linear regression}
\label{tb.effs}
\end{center}
\end{table}

\begin{remark}
In simple linear regression models, the works by Neter and Montgomery \cite{mongto, neter} considered the correlation of $\hat{\theta_{1}}$ and $\hat{\theta_{2}}$ in the context of simultaneous inferences about these two parameters, where they highlight that the correlation of the estimators leads to over or underestimation of the parameters as a function of the elongation of the confidence ellipse. 
In summary, it can be said that the discussion about the most efficient way to simultaneously control precision and correlation is not finished. It is necessary to revisit these concepts and evaluate new alternatives that alone or in synergy with other selection criteria allow researchers to make optimal decisions.
\end{remark}

\section{Case study: Michaelis-Menten Model}
\subsection{Correlation of parameter estimators}\label{parcor}
The use of the Michaelis-Menten model has covered a great diversity of disciplines, chemistry, ecology, epidemiology, among many. Where the development and theoretical interpretation of the parameters allows to explain the behavior of the phenomenon represented in its two parameters \textit{V} and \textit{K}, 
\begin{equation}
\label{eq:mmm}
E\left[ v \right]=\frac{Vx}{K+x},\; \text{var}(v)=\sigma ^{2},\, x\in\chi =\left[ 0,bK \right].
\end{equation}
Usually $v$ represents a rate and $x$ the variable to be controlled and it is assumed normally distributed. Parameter $V$ would be the maximum rate the system can reach asymptotically. Parameter $K$ is known as the Michaelis--Menten constant, so that when $x$ reaches the value $K$, the variable $v$ reaches half of its maximum $V$. Intuitively the values of $K$ and $V$ can be expected to present a high correlation. The design space here depends on the parameter $K$, which is quite realistic and frequent in the literature. On the other hand $K$ is unknown and it has to be estimated from the data. For the computation of an optimal experimental design a nominal value will be used for this and also for computing the FIM.

\subsection{Michaelis Menten correlation evidence}
Research papers that report problems of high correlation between the estimators of the parameters, while optimizing a  criterion for this model, are scarce. López-Fidalgo \& Wong (1998) \cite{fidalgo1} provided D-- and c--optimal designs for this model, as well as regular optimal sequences of design conditions. 

\begin{figure}
\centering
\includegraphics[trim={1mm 1mm 1mm 1mm},clip,width=12cm]{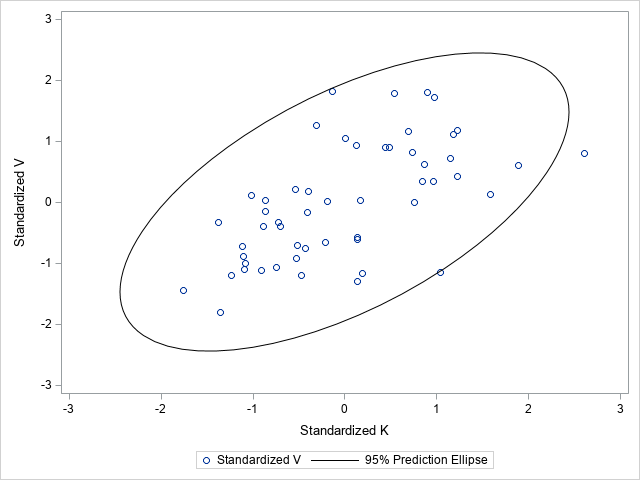}
\caption{Scatter plot of the standardized estimations of parameters $K$ and $V$ in six independent studies}
\label{f05ckv}
\end{figure}

Here, estimated parameter values for this model were taken from six independent studies in different research areas: 
plant physiology \cite{parmmfv}, cancer detection \cite{parmmcd}, phytoremediation \cite{parmmfr}, asthma \cite{parmmas}, and enzyme kinetic \cite{parmmki}. Then, the parameter estimators were standardized by subtracting their mean and dividing it by their standard deviation. This procedure was carried out separately for each parameter and for each study. Table \ref{t01corkv} shows the number of datasets in each study and the correlation between the estimated parameters. A general correlation of $r= 0.69$ was obtained for all the estimators in all the studies. Figure  \ref{f05ckv} shows the standardized values of the estimators $\hat{V}$ versus $\hat{K}$ with the empirical confidence ellipse. 
All the studies led to high correlations, but not enough to create invertibility problems of the FIM. However, these values warn of a real situation that was not considered in studies on optimal designs for this model and that may affect the best performance of the recommended proposals. 

\begin{table}[H]
 \begin{center}
 \begin{tabular}{l r l } \hline
 Study & $N^a$ & $r$ \\ \hline
 Hasegawa y Ichii, 1994 & 16 & 0.70 \\ 
 Blokh et al, 2007  & 9 & 0.88 \\
 Yu et al, 2005 (Lineweaver-Burk method)  & 12 & 0.52 \\
 Ishizaki, Kubo 1987(a) & 6 & 0.81 \\
 Ishizaki, Kubo 1997(b) & 5 & 0.89 \\
 Goudar et al, 1999  & 6 & 0.81 \\ \hline
 \end{tabular}\\
 \footnotesize{$^a$ number of different values estimated of $K$ and $V$}
 \caption{Empirical evidence of correlation in parameter estimators for the Michaelis-Menten model for data from five  different papers \label{t01corkv}}
 \end{center}
\end{table}


\subsection{Criteria for simultaneous control of correlation and estimation precision}
The performance of R--optimality compounded with some other criteria will be to evaluated for the Michaelis Menten model. Compound criteria that have been previously used by other authors for analyzing their behavior in relation to the level of correlation of the estimators of the parameters will be used. In particular:
\begin{enumerate}
\item A-optimality is an appropriate candidate to be considered jointly with correlation since it only considers the trace of the inverse of the FIM. To avoid biases by orders of magnitude between the parameter estimators, the standardized version, SA--optimality, will be used, 

\begin{equation}
\label{19_aop}
\Phi_{SA}^{*}(\xi) = \frac{\Phi_{c_{1}}(\xi)}{\Phi_{c_{1}}(\xi_{c_{1}}^{*})}+\frac{\Phi_{c_{2}}(\xi)}{\Phi_{c_{2}}(\xi_{c_{2}}^{*})},
\end{equation}
where

\begin{equation}
\Phi_{c}(\xi) =c^{t}M^{-}(\xi)c
\end{equation}
stands for the generalized variance with $c\in R^{2}$ and $c_1^T=(1,0)$, $c_2^T=(0,1)$. The superscript $-$ stands for the class of pseudo--inverses of a matrix. The linear combination of parameters $c^T\theta$ is estimable if, and only if $c^{t}M^{-}(\xi)c$ is constant for any member of the class of pseudo--inverses. Design $\xi_c^\star$ is a c--optimal design. Notice that

\begin{equation*}
\begin{matrix}
\text{if } c=c_{1}=\begin{pmatrix}
1\\ 
0
\end{pmatrix}\; \text{then} \; \Phi_{c_{1}}(\xi) =\left \{ M^{-}(\xi) \right \}_{11},\\ 
\text{if } c=c_{2}=\begin{pmatrix}
0\\ 
1
\end{pmatrix}\; \text{then} \; \Phi_{c_{2}}(\xi) =\left \{ M^{-}(\xi) \right \}_{22}.
\end{matrix}
\end{equation*}

\item D-optimality was already used by \cite{duggleby, fidalgo1} for the Michaelis-Menten model. The D--optimal design is always,

\begin{equation}
\label{dimm}
\xi _{D}=
\begin{Bmatrix}
\frac{b}{2+b}K &bK \\ 
1/2 & 1/2
\end{Bmatrix}
\end{equation}
Nominal values of the parameters $V=43.73$ and $K=227.27$ are taken from \cite{duggleby, fidalgo1} for later efficiency calculations.
\item E-optimality has been one of the criteria specifically considered to deal with the problem of the correlation of the parameter estimators. This criterion seeks to minimize the length of the major axis of the confidence ellipsoid of the parameters. However, instead of taking into account the major axis, the ratio between the major and the minor axes lengths, the so called condition number, should be more efficient for optimality purposes (E modified or EM--optimality).
\item A general version of $r^{2}$-optimality for more than two parameters 
had already been provided by Pritchard \& Bacon \cite{pritchard_bacon}, who proposed a correlation-specific optimality criterion (\ref{06_cpb}), which measures the size of the correlation. They called it C--optimality,
\begin{equation}
\label{06_cpb}
\Phi_{c}=\left \{ \frac{\sum_{ij; i\neq j}r_{ij}^{2}}{p(p-1)} \right \}^{1/2},
\end{equation}
\end{enumerate}
where $r_{ij}$ is the correlation coefficient between estimators of the parameters $i$ and $j$ and $p$ the number of parameters.

The design space considered by Dugglevy \cite{duggleby} was $[0,bK]$ with values of $b$ between $0.25$ and $4$. Given the typical convergence problems near $x=0$, the search for the optimal design was restricted to design spaces $[\epsilon,bK]$ with $\epsilon=0.05, 0.5, 1$. Table~\ref{t02dis} shows how convergence problems near the origin ($\epsilon=0$) lead to single--point designs for the criteria $\phi_{EM}$ and $\phi_{r^{2}}$, contrary to the criteria $\phi_{D}$, $\phi_{SA}$ and $\phi_{R}$, which obtained very similar results with two-point designs. As the value of $\epsilon$ increases, all designs tend to identify models with more than one point, with the exception of EM--optimality, which gives two-point designs for $\epsilon$=0.5 and 1, although with weights concentrated at the initial design point. 
For these computations $b=5$ was chosen according to some literature.

\begin{table}
\begin{center}
\begin{tabular}{llll}
\\ \hline
$\epsilon$*&Criteria&a&p\\ \hline
$\epsilon$=0&$\phi_{D}$&0.71 &0.50\\
&$\phi_{SA}$&0.55 &0.54\\
&$\phi_{R}$&0.55 &0.54\\
&$\phi_{EM}$&0.00 &1.00\\
&$\phi_{r^{2}}$&0.00 &1.00\\
$\epsilon$=0.05&$\phi_{D}$&0.71 &0.50\\
&$\phi_{SA}$&0.55 &0.54\\
&$\phi_{R}$&0.55 &0.53\\
&$\phi_{EM}$&0.05 &1.00\\
&$\phi_{r^{2}}$&0.50 &0.61\\
$\epsilon$=0.5&$\phi_{D}$&0.71 &0.50\\
&$\phi_{SA}$&0.55 &0.54\\
&$\phi_{R}$&0.55 &0.53\\
&$\phi_{EM}$&0.50 &0.86\\
&$\phi_{r^{2}}$&0.50 &0.61\\
$\epsilon$=1&$\phi_{D}$&1.00 &0.50\\
&$\phi_{SA}$&1.00 &0.50\\
&$\phi_{R}$&1.00 &0.49\\
&$\phi_{EM}$&1.00 &0.73\\
&$\phi_{r^{2}}$&1.00 &0.48\\ \hline
 \end{tabular}
 \caption{Optimal designs for estimating the parameters $K$ and $V$,
for different criteria and values of $\epsilon$. The support points of these designs are $aK$ and $5K$ with mass $p$ at $aK$}
\label{t02dis}
\end{center}
\end{table}
Table~\ref{t03efi} provides cross--efficiencies and the square of the correlations for the criteria considered. It is noted that none of the criteria (even $\phi_{r^{2}}$) manage to reduce the correlation to zero. The criteria that achieve the best reduction in correlation are $\phi_{R}$ and, unexpectedly, $\phi_{SA}$. It is also found that the increase in the value of $\epsilon$ causes the correlation levels to increase, although they are relatively stable when they oscillate between 0.05 and 0.5. Regarding the efficiencies, there are no big differences between the different designs. The EM--optimal design is the one with the worst performance at any value of $\epsilon$.

\begin{table}
\begin{center}
\begin{tabular}{lllllll|l}
\\ \hline
$\epsilon$&Criteria& $\text{Eff}_{D}$ & $\text{Eff}_{SA}$ & $\text{Eff}_{R}$ & $\text{Eff}_{EM}$ & $\text{Eff}_{r^{2}}$ & $r^{2}$\\ \hline
0&$\phi_{D}$&1.00&0.99&0.99& -- &0.00&0.69\\
&$\phi_{SA}$&0.97&1.00&1.00& -- &0.00&0.66\\
&$\phi_{R}$&0.97&1.00&1.00& -- &0.00&0.66\\
&$\phi_{EM}$& -- & -- & -- &1.00& -- & -- \\
&$\phi_{r^{2}}$&0.71&0.75&0.75& -- &1.00&0.49\\
0.05&$\phi_{D}$&1.00&0.99&0.99&0.16&0,97&0,69\\
&$\phi_{SA}$&0.97&1.00&1.00&0.02&0.98&0.66\\
&$\phi_{R}$&0.97&1.00&1.00&0.03&0.98&0.66\\
&$\phi_{EM}$&0.18&0.22&0.22&1.00&0.26&0.66\\
&$\phi_{r^{2}}$&0.92&0.98&0.98&0.96&1.00&0.64\\
0.5&$\phi_{D}$&1.00&0.99&0.99&0.82&0.97&0.69\\
&$\phi_{SA}$&0.97&1.00&1.00&0.55&0.98&0.66\\
&$\phi_{R}$&0.97&1.00&1.00&0.60&0.98&0.66\\
&$\phi_{EM}$&0.45&0.53&0.53&1.00&0.64&0.73\\
&$\phi_{r^{2}}$&0.92&0.98&0.98&0.87&1.00&0.64\\
1&$\phi_{D}$&1.00&1.00&1.00&0.94&1.00&0.75\\
&$\phi_{SA}$&1.00&1.00&1.00&0.78&1.00&0.75\\
&$\phi_{R}$&1.00&1.00&1.00&0.79&1.00&0.75\\
&$\phi_{EM}$&0.78&0.78&0.77&1.00&0.75&0.80\\
&$\phi_{r^{2}}$&1.00&1.00&1.00&0.94&1.00&0.75\\ \hline
 \end{tabular}
 \caption{Efficiencies and correlations for the criteria evaluated with different values of $\epsilon$}
 \label{t03efi}
 \end{center}
\end{table}

\subsection{Compound optimal criteria and Pareto front for controlling the correlation and estimation precision}
Based on the previous results, it was decided to consider compound criteria. The results of EM-- and D--optimal designs was unexpected showing the worst control on the correlation of the estimators of the parameters. Surprisingly A--optimal designs gave the lowest correlation values. This, together with the problems related to singular EM--optimal designs led to the consideration of D--optimality as the basis for compound criteria with R--optimality,
\begin{equation}
\label{08_cef}
\Phi_{\lambda } \left (\xi \right )= \left ( 1-\lambda \right ) \frac{1}{\text{Eff}_{D}(\xi)}+ \lambda\,\frac{1}{\text{Eff}_{R}(\xi)} \; \; \; \: 
\lambda\; \epsilon\left [ 0,1 \right ],
\end{equation}
where $0\leq \lambda\leq1$ is a user-selected constant \cite{fidalgo1}. Since $\phi_{D}$, $\phi_{SA}$ and $\phi_{R}$ are convex criteria, therefore, linear combinations of them are still convex, and the equivalence theorem is valid. Additionally, they are also inverse positive homogeneous.

\begin{figure}
\centering
\includegraphics[trim={1mm 1mm 1mm 1mm},clip,width=7cm]{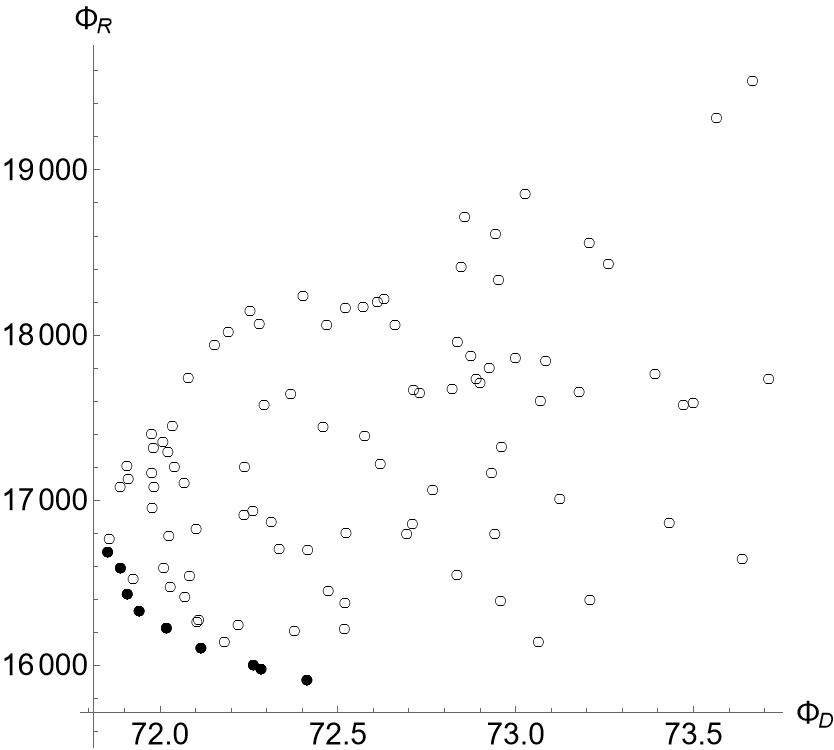}
\caption{Pareto front for $\Phi_{R}$ y $\Phi_{D}$ criteria}
\label{f06par}
\end{figure}

\begin{figure}
\centering
\includegraphics[trim={1mm 1mm 1mm 1mm},clip,width=7cm]{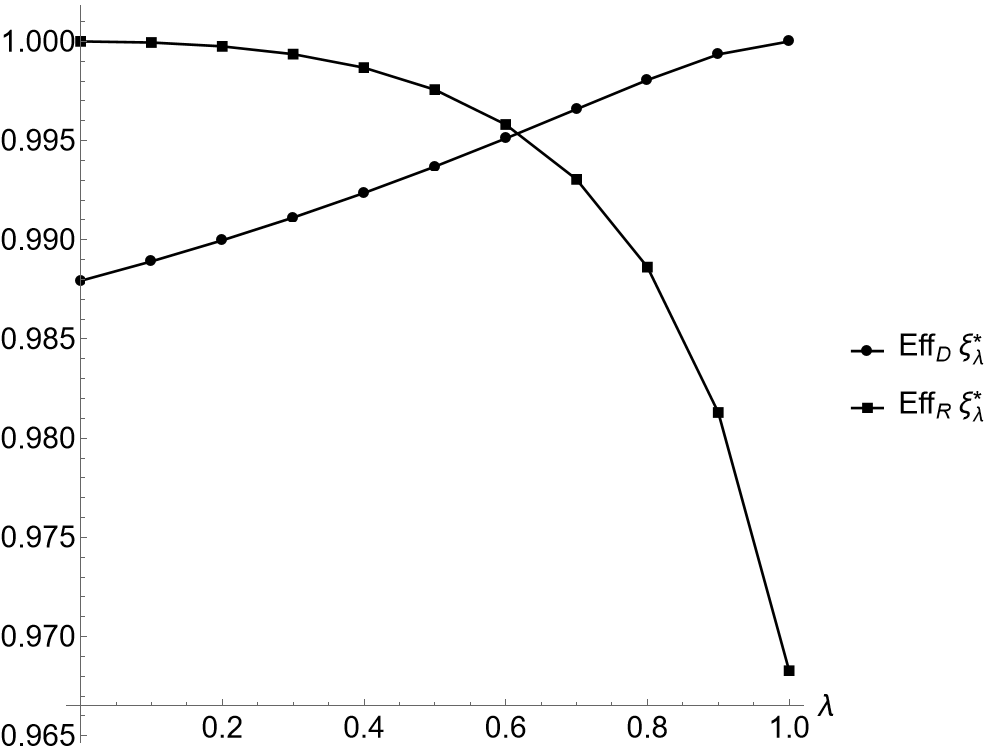}
\caption{Compound optimal design}
\label{f07efird}
\end{figure}

The Pareto front was built with 1000 simulated two--point designs. From them nine admissible designs were obtained (Figure \ref{f06par}). D--efficiency and correlation of these designs is shown in Table \ref{t04efipar}. Then the efficiencies of the selected designs were plotted using the compound criterion of $\Phi_{R}$ and $\Phi_{D}$. The lowest correlation design was obtained from the Pareto front 
(Figure \ref{f07efird}).

\begin{table}
 \begin{center}
 \resizebox{15cm}{!}{
 \begin{tabular}{cccccccccccc}\hline
 & $\xi_{D}^{*}$ & $\xi_{R}^{*}$ & $\xi_{DR_{1}}^{*}$ & $\xi_{DR_{2}}^{*}$ & $\xi_{DR_{3}}^{*}$ & $\xi_{DR_{4}}^{*}$ & $\xi_{DR_{5}}^{*}$ & $\xi_{DR_{6}}^{*}$ & $\xi_{DR_{7}}^{*}$ & $\xi_{DR_{8}}^{*}$ & $\xi_{DR_{9}}^{*}$ \\ \hline
 $\text{Eff}_{D}$ & 1 & 0.968 & $ 1.000 $ & $1.000$ & $1.000$ & $1.000$ & $1.000$ & $1.000$ & $1.000$ & $1.000$ & $1.000$ \\
 \text{p} & 0.50 & 0.465 & 0.492 & 0.485 & 0.480 & 0.494 & 0.482 & 0.482 & 0.501 & 0.494 & 0.489 \\
 \text{a} & 0.714 & 0.551 & 0.706 & 0.699 & 0.695 & 0.709 & 0.696 & 0.697 & 0.715 & 0.708 & 0.703 \\
 $r^{2}$ & 0.692 & 0.656 & 0.690 & 0.688 & 0.687 & 0.691 & 0.688 & 0.688 & 0.693 & 0.691 & 0.689 \\ \hline
 \end{tabular}
 }
 \caption{Designs selected from Pareto front for $\epsilon=0.5
$\label{t04efipar}}
 \end{center}
\end{table}

\subsection{Functional relationship between criteria}
\begin{figure}
\begin{center}
\subfigure{\includegraphics[width=55mm]{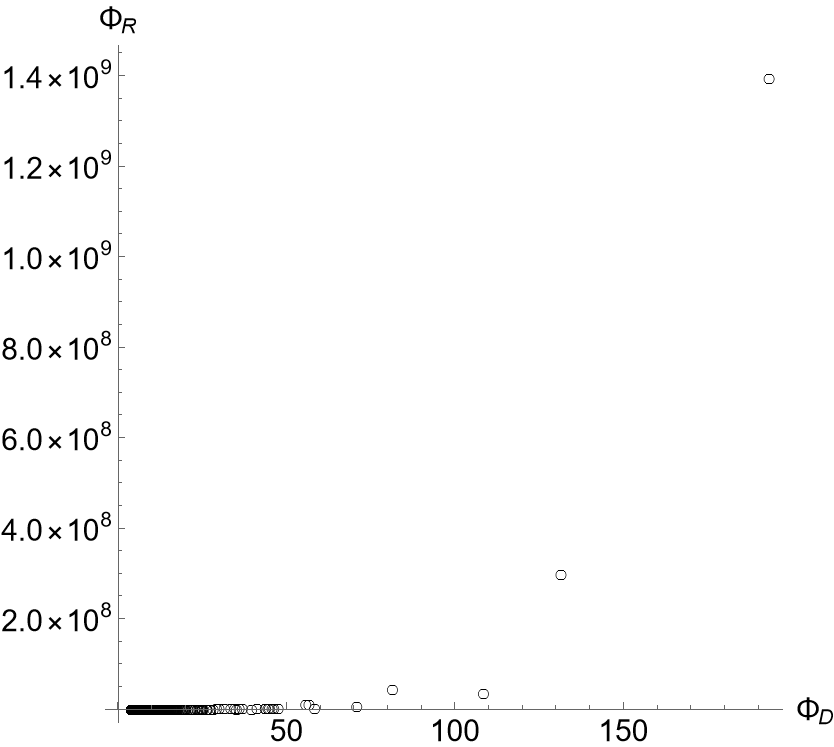}}
\subfigure{\includegraphics[width=55mm]{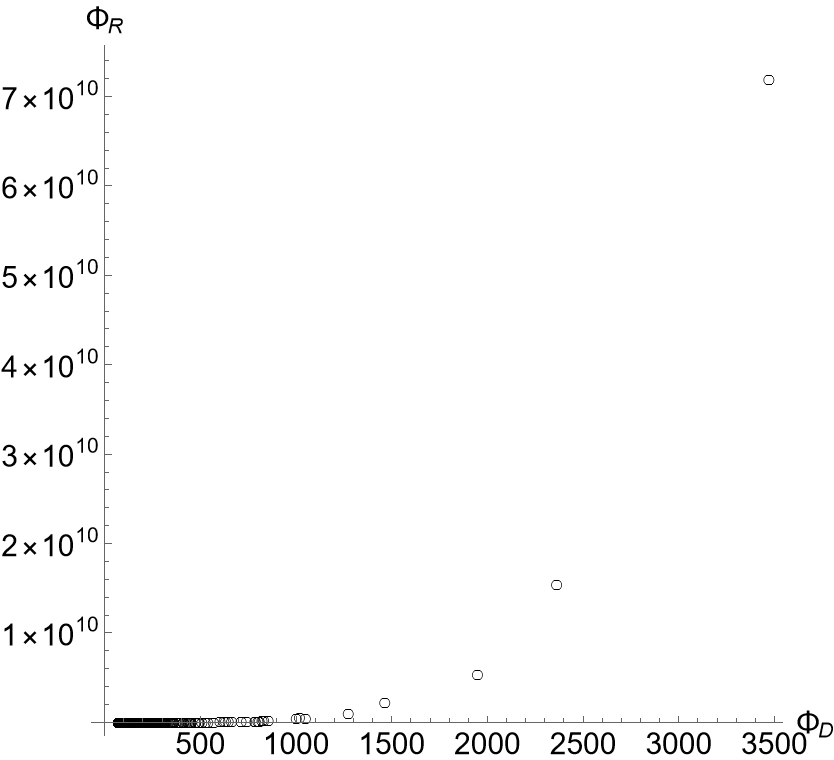}}
\subfigure{\includegraphics[width=55mm]{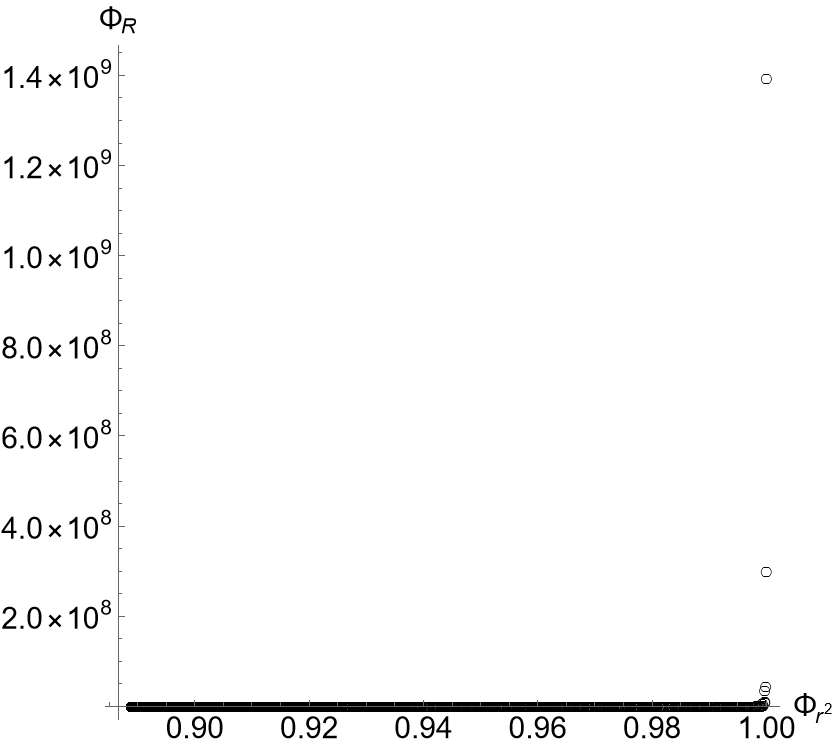}}
\subfigure{\includegraphics[width=55mm]{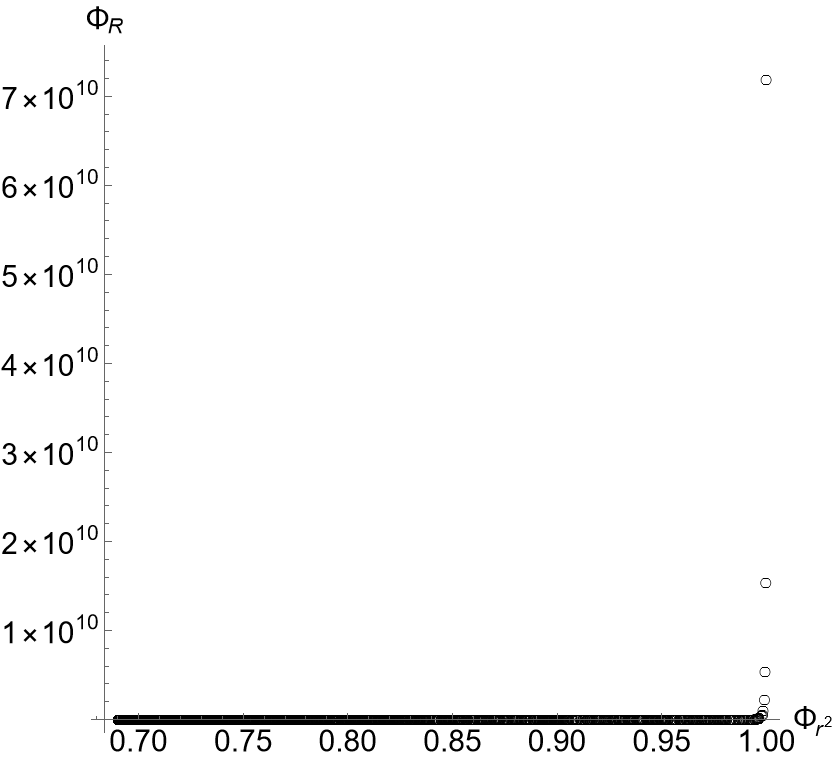}}
\subfigure{\includegraphics[width=55mm]{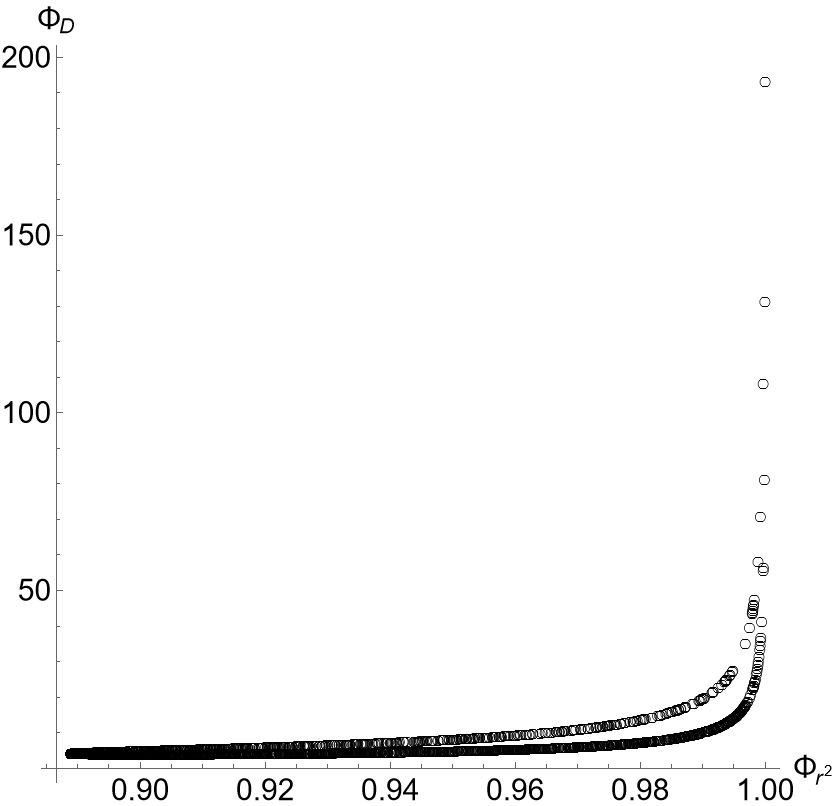}}
\subfigure{\includegraphics[width=55mm]{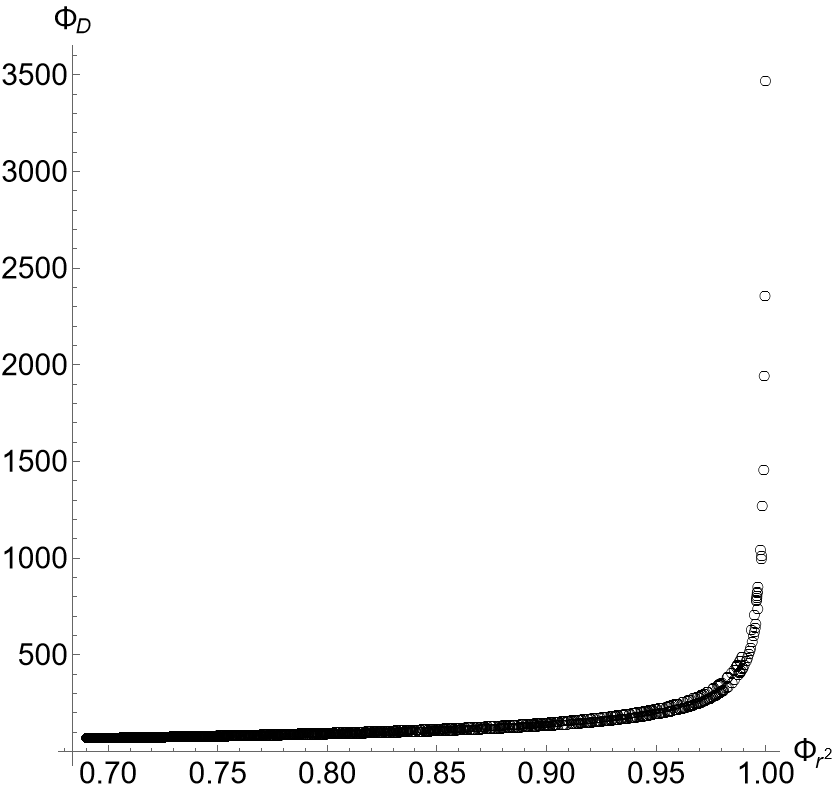}}
\caption{Functional relationship between values of $\phi_{R}$ versus $\phi_{R}$ and $\phi_{r^{2}}$ for simple linear regression fixing the first point at $a=0.5$ (left) and for the Michaelis Menten model fixing the first point at $a=0.71$ (right) varying $p\in (0,1)$ \label{f08rfc}}
\end{center}
\end{figure}

\begin{figure}
\begin{center}
\subfigure{\includegraphics[width=55mm]{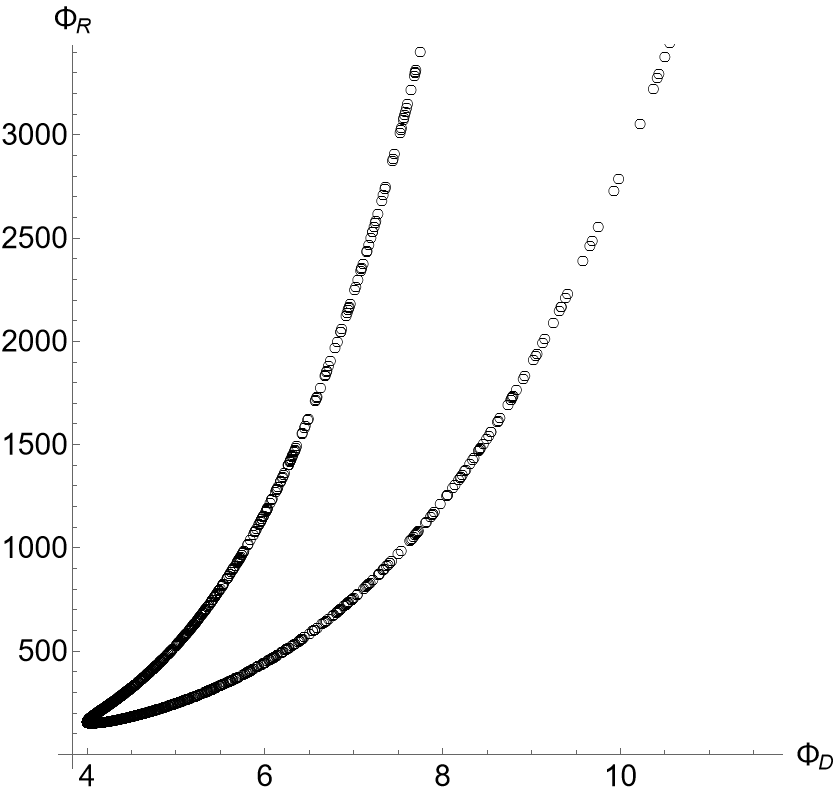}}
\subfigure{\includegraphics[width=55mm]{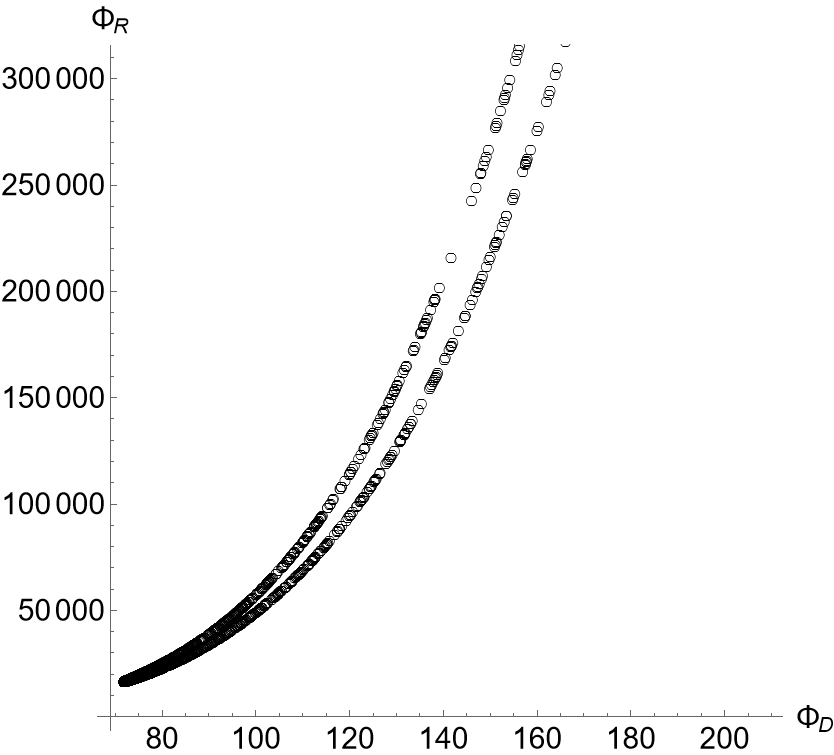}}
\subfigure{\includegraphics[width=55mm]{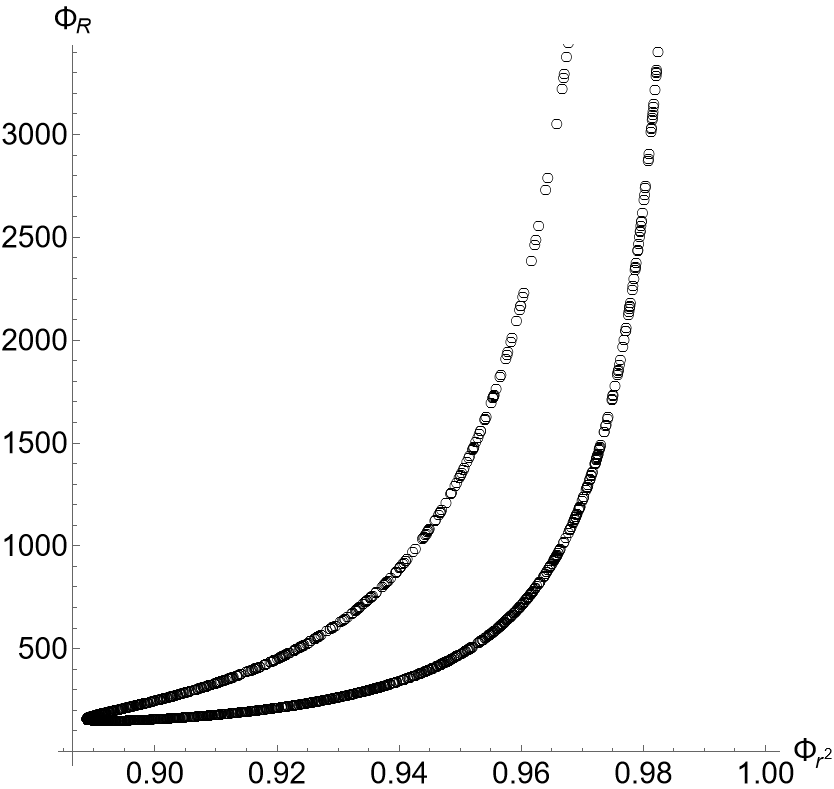}}
\subfigure{\includegraphics[width=55mm]{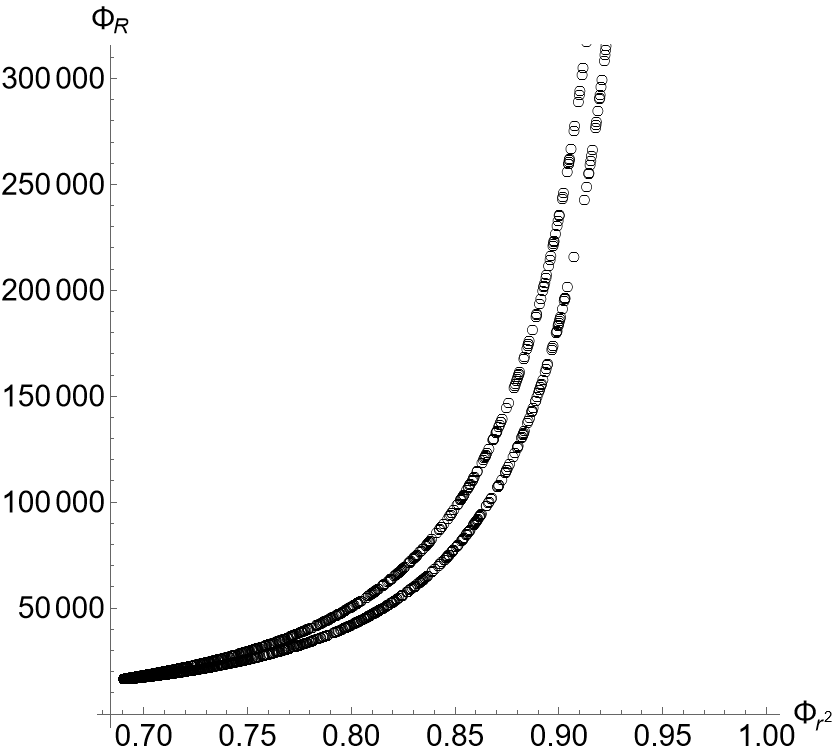}}
\subfigure{\includegraphics[width=55mm]{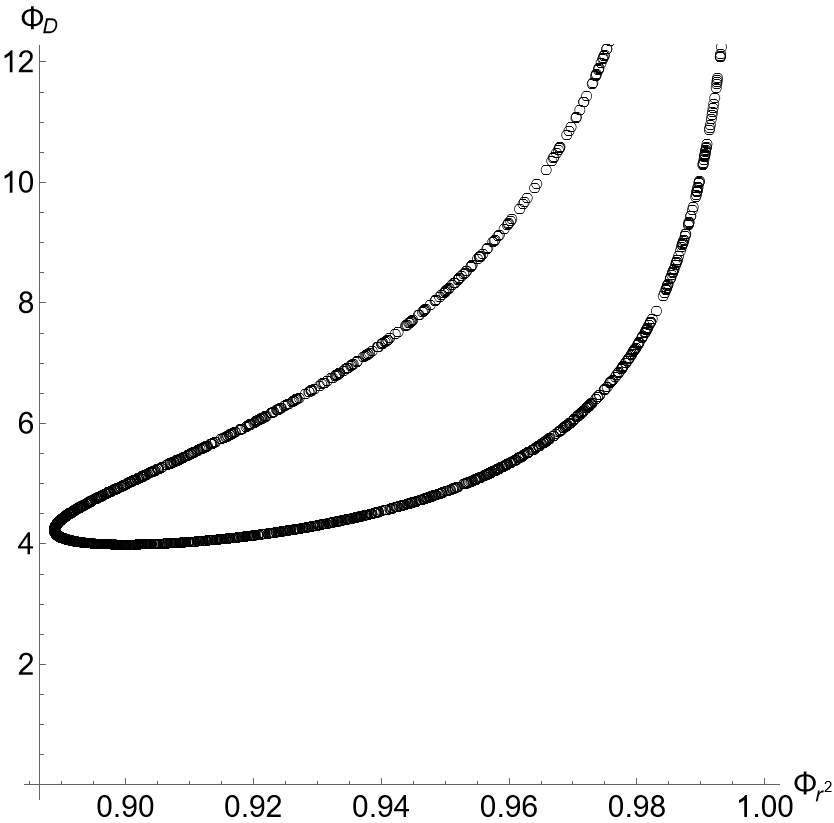}} 
\subfigure{\includegraphics[width=55mm]{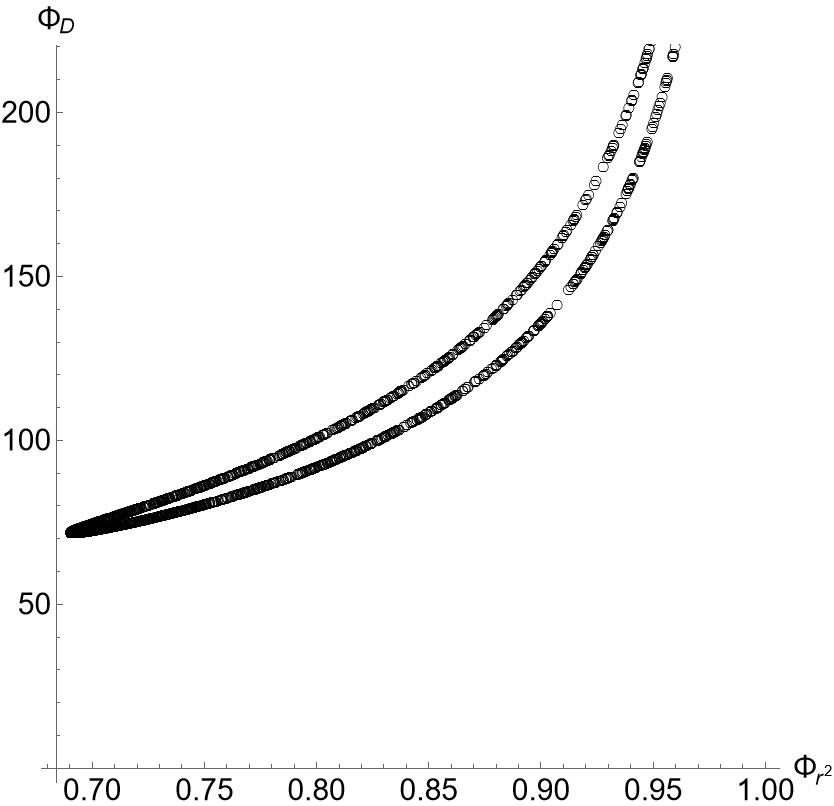}}
\caption{Functional relationship between values of $\phi_{R}$ versus $\phi_{R}$ and $\phi_{r^{2}}$ for simple linear regression fixing the first point at $a=0.5$ (left) and for the Michaelis Menten model fixing the first point at $a=0.71$ (right) varying $p\in (0,1)$ with a magnifying glass effect near the optimal designs 
\label{f09rfc}}
\end{center}
\end{figure}

\begin{figure}
\subfigure[Simple linear regression model]{\includegraphics[width=70mm]
{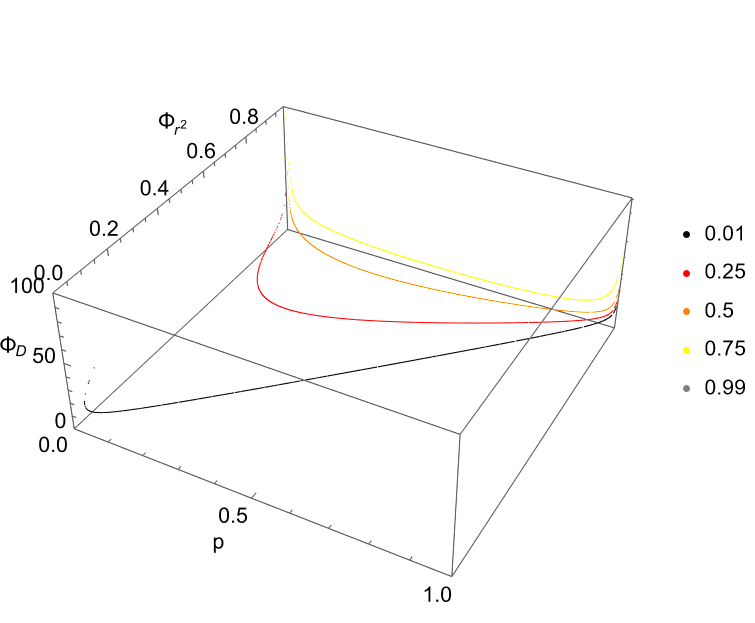}}
\subfigure[Michaelis Menten model]{\includegraphics[width=70mm]{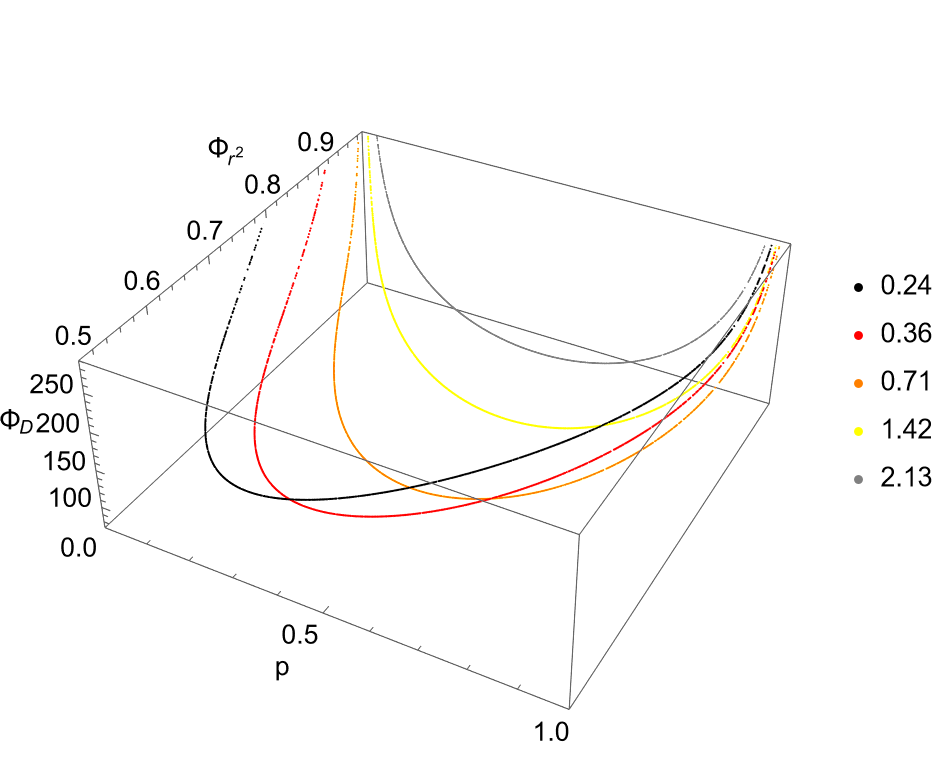}}
\caption{Loop effect due to imbalance in the masses of the D--optimal designs for simple linear regression from $a=0.01$ to $0.99$; and for the Michaelis Menten model from $a=0.24$ to $2.13$, and varying $p\in (0,1)$ with a magnifying glass effect near the optimal designs}
\label{f010loop}
\end{figure}

To evaluate the functional relationship between the three criteria $\phi_{R}$, $\phi_{r^{2}}$ and $\phi_{D}$ in Equation (\ref{eq:03_cdr}), two--point designs with a weight $p\in (0,1)$ and fixing $a=0.71$ for the Michaelis-Menten model and $a=0.5$ for the linear model, were simulated. It can be seen that the rate of change of $\phi_{D}$ and $\phi_{R}$, in relation to $\phi_{r^{2}}$, both for the Michaelis Menten and for the linear models, is different. It is observed that for high values of correlation ($r^{2}>$0.99) $\phi_{R}$ converges faster to a minimum than $\phi_{D}$ (Figure \ref{f08rfc}). On the other hand in Figure \ref{f09rfc}, where the correlation converges to minimum values, the expected behavior is observed in the form of a loop (Figure \ref{f010loop}). This corresponds to the Pareto front principle: there are designs for which each criterion is not minimal with no better solutions for the other criterion.

\section{Discussion}\label{discussion}
The correlation of the estimators of the parameters is, in addition to the precision of estimation, one of the necessary conditions to look for optimal designs of experiments. In this way, its consideration has always been implicit in the definition of any design, although, it is necessary to clarify that it has always been of low priority in relation to the criteria related to the improvement in the precision of the estimation.

In the present work, previous strategies in the literature for this purpose were reviewed. They can be classified into three groups: 1) those considering that the alphabetic criteria look at the same time for the efficiency in the estimation of the parameters and control the correlation of their estimators; 2) those that seek for criteria focused on spherical reliable confidence ellipsoids; and 3) those using the combination of different criteria to balance between these two characteristics.

In this way, a new criterion was found that corresponds to the third approach (simultaneous control of the efficiency in the estimation of the parameters and reduction in the correlation of the estimators). This new criterion coincides with the proposal given by Dette (1997)\cite{dette_1}, and it has similar results than SA--optimality since it considers only the information from the diagonal of the inverse of the information matrix, which in principle would be a limitation for its use in optimization of the correlation.

The evaluations carried out on a Michaelis Menten model showed that criteria such as D-- or EM--optimality, which have been reported to control the correlation of the estimators of the parameters, showed a poor performance. On the contrary, SA-- and R--optimality criteria provided lower correlation values.

The first comparison analysis between the different criteria evaluated was the computation of different optimal designs, from which it can be observed that if the experimental area includes zero, the optimal designs for the $\phi_{EM}$ and $\phi_{r^{2}}$ criteria are singular with maximum weights at the origin. The D--optimal design is not affected by the evaluated region and SA-- and R--optimal designs gave slightly lower correlation than the D--optimal design, but with greater weight towards the initial point. Conversely, if the value of lowest possible value is set to a value less or equal to 1, the same points are identified for the different criteria. With the exception of the EM--optimal design, the weight of the design tends to the lower extreme point. This has a similar behavior for the case with the lower extreme of the design space equal to zero. At the intermediate values  of this lower extreme from 0.05 to 0.5, both the locations and the weights are similar between the different designs, except for the EM--optimal design which generally tends to concentrate the weights at the lower extreme point. From this analysis it can be seen that, for practical reasons, the value of 0.5 for the lower extreme point can be considered a sufficient exploration limit for the design space. The efficiency calculations confirm what has been found in the definition of optimal designs: the extreme values (0 and 1) are not good optimality evaluation points, while the intermediate value of 0.05 gives the best efficiency values for all criteria.

The results showing the lower correlation values are obtained for SA-- and R--optimal design while the worst are for the EM--optimal design seems to contradict the fact that SA-- and R--optimality only take into account the diagonal of the inverse of the information matrix and not the covariance. Considering the geometric properties of the confidence ellipsoid, from which all the criteria evaluated are related, the product of the variances ($\phi_{R}$ criterion) seems to recover the information on the correlation between the estimators better than other criteria.

Multi--objective optimization either using \textit{Compound optimal design} or using Pareto fronts tend to give more weight to estimation efficiency than to correlation, or to give solution ranges where the final selection becomes a subjective decision of the researcher. In this sense, R--optimality frees the researcher from that decision, selecting a design with high efficiency values for other criteria, such as D--optimality, with the lowest possible level of correlation.

As mentioned in the Introduction some authors  have claimed for many years that minimizing the determinant of the inverse of the FIM simultaneously reduces the level of correlation as well as optimizing the precision of the estimators. The truth is that this only happens when the designs are far from the D--optimal design. Nevertheless, close to the D--optimal design, the relationship between correlation and precision show a "loop effect" (Figure \ref{f010loop}), which allows the Pareto front to generate different solutions. 

\section*{Acknowledgement}
This work was sponsored by Ministerio de Ciencia e Innovación, PID2020-113443RB-C21.

\bibliographystyle{elsarticle-harv} 
\bibliography{cas-refs}

\begin{thebibliography}{25}
\expandafter\ifx\csname natexlab\endcsname\relax\def\natexlab#1{#1}\fi
\providecommand{\url}[1]{\texttt{#1}}
\providecommand{\href}[2]{#2}
\providecommand{\path}[1]{#1}
\providecommand{\DOIprefix}{doi:}
\providecommand{\ArXivprefix}{arXiv:}
\providecommand{\URLprefix}{URL: }
\providecommand{\Pubmedprefix}{pmid:}
\providecommand{\doi}[1]{\href{http://dx.doi.org/#1}{\path{#1}}}
\providecommand{\Pubmed}[1]{\href{pmid:#1}{\path{#1}}}
\providecommand{\bibinfo}[2]{#2}
\ifx\xfnm\relax \def\xfnm[#1]{\unskip,\space#1}\fi
\bibitem[{Agarwal and Brisk(1985)}]{agarwal}
\bibinfo{author}{Agarwal, A.}, \bibinfo{author}{Brisk, M.},
  \bibinfo{year}{1985}.
\newblock \bibinfo{title}{Sequential experimental design for precise parameter
  estimation. 1. use of reparameterization}.
\newblock \bibinfo{journal}{Industrial \& Engineering Chemistry Process Design
  and Development} \bibinfo{volume}{24}, \bibinfo{pages}{203--207}.
\bibitem[{Bhonsale et~al.(2022)Bhonsale, Nimmegeers, Akkermans, Telen, Stamati,
  Logist and Van-Impe}]{bhonsale}
\bibinfo{author}{Bhonsale, S.}, \bibinfo{author}{Nimmegeers, P.},
  \bibinfo{author}{Akkermans, S.}, \bibinfo{author}{Telen, D.},
  \bibinfo{author}{Stamati, I.}, \bibinfo{author}{Logist, F.},
  \bibinfo{author}{Van-Impe, J.}, \bibinfo{year}{2022}.
\newblock \bibinfo{title}{Optimal experiment design for dynamic processes. In
  Simulation and Optimization in Process Engineering}.
\newblock \bibinfo{publisher}{Elsevier}.
\bibitem[{Blokh et~al.(2007)Blokh, Stambler, Afrimzon, Shafran, Korech,
  Sandbank, Orda, Zurgil and Deutsch}]{parmmcd}
\bibinfo{author}{Blokh, D.}, \bibinfo{author}{Stambler, I.},
  \bibinfo{author}{Afrimzon, E.}, \bibinfo{author}{Shafran, Y.},
  \bibinfo{author}{Korech, E.}, \bibinfo{author}{Sandbank, J.},
  \bibinfo{author}{Orda, R.}, \bibinfo{author}{Zurgil, N.},
  \bibinfo{author}{Deutsch, M.}, \bibinfo{year}{2007}.
\newblock \bibinfo{title}{The information-theory analysis of michaelis–menten
  constants for detection of breast cancer}.
\newblock \bibinfo{journal}{Cancer Detection and Prevention}
  \bibinfo{volume}{31}, \bibinfo{pages}{489--498}.
\bibitem[{Box and Hunter(1963)}]{boxhunter}
\bibinfo{author}{Box, G.}, \bibinfo{author}{Hunter, W.}, \bibinfo{year}{1963}.
\newblock \bibinfo{title}{Sequential design of experiments for non-linear
  models}, in: \bibinfo{booktitle}{Proceedings of the IBM Scientific Computing
  Symposium on Statistics}, pp. \bibinfo{pages}{113--137}.
\bibitem[{Das and Dennis(1998)}]{das}
\bibinfo{author}{Das, I.}, \bibinfo{author}{Dennis, J.}, \bibinfo{year}{1998}.
\newblock \bibinfo{title}{Normal-boundary intersection: A new method for
  generating the pareto surface in nonlinear multicriteria optimization
  problems}.
\newblock \bibinfo{journal}{SIAM journal on optimization} \bibinfo{volume}{8},
  \bibinfo{pages}{631--657}.
\bibitem[{Dette(1997)}]{dette_1}
\bibinfo{author}{Dette, H.}, \bibinfo{year}{1997}.
\newblock \bibinfo{title}{Designing experiments with respect to
  ‘standardized’ optimality criteria}.
\newblock \bibinfo{journal}{Journal of the Royal Statistical Society: Series B
  (Statistical Methodology)} \bibinfo{volume}{59}, \bibinfo{pages}{97--110}.
\bibitem[{Duggleby(1979)}]{duggleby}
\bibinfo{author}{Duggleby, R.}, \bibinfo{year}{1979}.
\newblock \bibinfo{title}{Experimental designs for estimating the kinetic
  parameters for enzyme-catalysed reactions}.
\newblock \bibinfo{journal}{Journal of Theoretical Biology}
  \bibinfo{volume}{81}, \bibinfo{pages}{671--684}.
\bibitem[{Franceschini and Macchietto(2008a)}]{stateart}
\bibinfo{author}{Franceschini, G.}, \bibinfo{author}{Macchietto, S.},
  \bibinfo{year}{2008}a.
\newblock \bibinfo{title}{Model-based design of experiments for parameter
  precision: State of the art}.
\newblock \bibinfo{journal}{Chemical Engineering Science} \bibinfo{volume}{63},
  \bibinfo{pages}{4846–4872}.
\bibitem[{Franceschini and Macchietto(2008b)}]{Franceschini}
\bibinfo{author}{Franceschini, G.}, \bibinfo{author}{Macchietto, S.},
  \bibinfo{year}{2008}b.
\newblock \bibinfo{title}{Novel anticorrelation criteria for model‐based
  experiment design: Theory and formulations}.
\newblock \bibinfo{journal}{AIChE Journal} \bibinfo{volume}{54},
  \bibinfo{pages}{1009–1024}.
\bibitem[{G. and H.(1959)}]{box_lucas}
\bibinfo{author}{G., B.}, \bibinfo{author}{H., L.}, \bibinfo{year}{1959}.
\newblock \bibinfo{title}{Design of experiments in non-linear situations}.
\newblock \bibinfo{journal}{Biometrika} \bibinfo{volume}{46},
  \bibinfo{pages}{77--90}.
\bibitem[{Goudar et~al.(1999)Goudar, Sonnad and Duggleby}]{parmmki}
\bibinfo{author}{Goudar, C.}, \bibinfo{author}{Sonnad, J.},
  \bibinfo{author}{Duggleby, R.}, \bibinfo{year}{1999}.
\newblock \bibinfo{title}{Parameter estimation using a direct solution of the
  integrated michaelis-menten equation}.
\newblock \bibinfo{journal}{Biochimica et Biophysica Acta}
  \bibinfo{volume}{1429(2)}, \bibinfo{pages}{377--383}.
\bibitem[{Hasegawa and Ichii(1994)}]{parmmfv}
\bibinfo{author}{Hasegawa, H.}, \bibinfo{author}{Ichii, M.},
  \bibinfo{year}{1994}.
\newblock \bibinfo{title}{Variation in michaelis-menten kinetic parameters for
  nitrate uptake by the young seedlings in rice (oryza sativa l.)}.
\newblock \bibinfo{journal}{Japanese Journal of Breeding}
  \bibinfo{volume}{44(4)}, \bibinfo{pages}{382--386}.
\bibitem[{Ishizaki and Kubo(1987)}]{parmmas}
\bibinfo{author}{Ishizaki, T.}, \bibinfo{author}{Kubo, M.},
  \bibinfo{year}{1987}.
\newblock \bibinfo{title}{Incidence of apparent michaelis-menten kinetic
  behavior of theophylline and its parameters (vmax and km) among asthmatic
  children and adults}.
\newblock \bibinfo{journal}{Therapeutic drug monitoring}
  \bibinfo{volume}{9(1)}, \bibinfo{pages}{11--20}.
\bibitem[{Lopez-Fidalgo and Wong(2002)}]{fidalgo1}
\bibinfo{author}{Lopez-Fidalgo, J.}, \bibinfo{author}{Wong, W.},
  \bibinfo{year}{2002}.
\newblock \bibinfo{title}{Design issues for the michaelis-menten model}.
\newblock \bibinfo{journal}{Journal of Theoretical Biology}
  \bibinfo{volume}{215}, \bibinfo{pages}{1--11}.
\bibitem[{Maheshwari et~al.(24)Maheshwari, Rangaiah and
  Samavedham}]{maheshwari}
\bibinfo{author}{Maheshwari, V.}, \bibinfo{author}{Rangaiah, G.},
  \bibinfo{author}{Samavedham, L.}, \bibinfo{year}{24}.
\newblock \bibinfo{title}{Multiobjective framework for model-based design of
  experiments to improve parameter precision and minimize parameter
  correlation}.
\newblock \bibinfo{journal}{Industrial \& Engineering Chemistry Research}
  \bibinfo{volume}{52}, \bibinfo{pages}{8289--8304}.
\bibitem[{McLean and McAuley(2012)}]{mclean}
\bibinfo{author}{McLean, K.A.P.}, \bibinfo{author}{McAuley, K.B.},
  \bibinfo{year}{2012}.
\newblock \bibinfo{title}{Mathematical modelling of chemical
  processes—obtaining the best model predictions and parameter estimates
  using identifiability and estimability procedures}.
\newblock \bibinfo{journal}{The Canadian Journal of Chemical Engineering}
  \bibinfo{volume}{90}, \bibinfo{pages}{351--366}.
\bibitem[{Montgomery and Peck(1992)}]{mongto}
\bibinfo{author}{Montgomery, D.}, \bibinfo{author}{Peck, E.},
  \bibinfo{year}{1992}.
\newblock \bibinfo{title}{Introduction to linear regresssion analysis}.
\newblock \bibinfo{publisher}{John Wiley \& Sons, Inc.}
\bibitem[{Neter et~al.(1990)Neter, Wasserman and Kutner}]{neter}
\bibinfo{author}{Neter, J.}, \bibinfo{author}{Wasserman, W.},
  \bibinfo{author}{Kutner, M.}, \bibinfo{year}{1990}.
\newblock \bibinfo{title}{Applied linear statistical models}.
\newblock \bibinfo{publisher}{Richard D. Irwin, Inc.}
\bibitem[{Pritchard and Bacon(1978)}]{pritchard_bacon}
\bibinfo{author}{Pritchard, D.}, \bibinfo{author}{Bacon, D.},
  \bibinfo{year}{1978}.
\newblock \bibinfo{title}{Prospects for reducing correlations among parameter
  estimates in kinetic models}.
\newblock \bibinfo{journal}{Chemical Engineering Science}
  \bibinfo{volume}{33(11)}, \bibinfo{pages}{1539--1543}.
\bibitem[{Rodriguez-Fernandez et~al.(2006)Rodriguez-Fernandez, Egea and
  Banga}]{valcor}
\bibinfo{author}{Rodriguez-Fernandez, M.}, \bibinfo{author}{Egea, J.},
  \bibinfo{author}{Banga, J.}, \bibinfo{year}{2006}.
\newblock \bibinfo{title}{Novel metaheuristic for parameter estimation in
  nonlinear dynamic biological systems}.
\newblock \bibinfo{journal}{BMC Bioinformatics} \bibinfo{volume}{7},
  \bibinfo{pages}{483}.
\bibitem[{Telen et~al.(2012)Telen, Logist, Van-Derlinden, Tack and
  Van-Impe}]{Telen}
\bibinfo{author}{Telen, D.}, \bibinfo{author}{Logist, F.},
  \bibinfo{author}{Van-Derlinden, E.}, \bibinfo{author}{Tack, I.},
  \bibinfo{author}{Van-Impe, J.}, \bibinfo{year}{2012}.
\newblock \bibinfo{title}{Optimal experiment design for dynamic bioprocesses: a
  multi-objective approach}.
\newblock \bibinfo{journal}{Chemical Engineering Science} \bibinfo{volume}{78},
  \bibinfo{pages}{82--97}.
\bibitem[{Vaghi et~al.(2020)Vaghi, Rodallec, Fanciullino, Ciccolini, Mochel,
  Mastri, Poignard, Ebos and Benzekry}]{gomred}
\bibinfo{author}{Vaghi, C.}, \bibinfo{author}{Rodallec, A.},
  \bibinfo{author}{Fanciullino, R.}, \bibinfo{author}{Ciccolini, J.},
  \bibinfo{author}{Mochel, J.}, \bibinfo{author}{Mastri, M.},
  \bibinfo{author}{Poignard, C.}, \bibinfo{author}{Ebos, J.},
  \bibinfo{author}{Benzekry, S.}, \bibinfo{year}{2020}.
\newblock \bibinfo{title}{Population modeling of tumor growth curves and the
  reduced gompertz model improve prediction of the age of experimental tumors}.
\newblock \bibinfo{journal}{PLoS computational biology}
  \bibinfo{volume}{16(2)}, \bibinfo{pages}{24}.
\bibitem[{Vaibhav et~al.(2013)Vaibhav, Gade and Lakshminarayanan}]{corr}
\bibinfo{author}{Vaibhav, M.}, \bibinfo{author}{Gade, P.},
  \bibinfo{author}{Lakshminarayanan, S.}, \bibinfo{year}{2013}.
\newblock \bibinfo{title}{Multiobjective framework for model-based design of
  experiments to improve parameter precision and minimize parameter
  correlation}.
\newblock \bibinfo{journal}{Industrial \& Engineering Chemistry Research}
  \bibinfo{volume}{52(24)}, \bibinfo{pages}{8289--8304}.
\bibitem[{Wang et~al.(2018)Wang, Yue and Yu}]{wang}
\bibinfo{author}{Wang, K.}, \bibinfo{author}{Yue, H.}, \bibinfo{author}{Yu,
  H.}, \bibinfo{year}{2018}.
\newblock \bibinfo{title}{Optimal input design for reduction of parameter
  correlations}, in: \bibinfo{booktitle}{In 2018 24th International Conference
  on Automation and Computing (ICAC)}, \bibinfo{address}{Newcastle upon Tyne,
  UK}. pp. \bibinfo{pages}{1--6}.
\bibitem[{Yu et~al.(2005)Yu, Zhou, Zhou and Liu}]{parmmfr}
\bibinfo{author}{Yu, X.}, \bibinfo{author}{Zhou, P.}, \bibinfo{author}{Zhou,
  X.}, \bibinfo{author}{Liu, Y.}, \bibinfo{year}{2005}.
\newblock \bibinfo{title}{Cyanide removal by chinese vegetation, quantification
  of the michaelis-menten kinetics}.
\newblock \bibinfo{journal}{Environmental Science and Pollution Research}
  \bibinfo{volume}{12(4)}, \bibinfo{pages}{221--226}.

\end{thebibliography}
\end{document}